\magnification\magstep1

\input amssym


\bigskip

\centerline{\bf Weakly compact approximation in Banach spaces}

\bigskip

\centerline{Edward Odell$^*$ and Hans-Olav Tylli$^{**}$}

\hbox{}\footnote{}{$^*$ Research supported by the NSF.}

\hbox{}\footnote{}{$^{**}$Research supported by the Academy of Finland
Project  \# 53893.}

\bigskip

{\leftskip 0.5 cm \rightskip 0.5 cm \noindent  {\bf Abstract.}
{\it  The Banach space $E$ has the weakly compact approximation property
(W.A.P. for short)  if there is a constant $C < \infty$
so that  for any weakly compact  set $D \subset E$ and $\varepsilon > 0$
there is a weakly compact operator $V: E \to E$ satisfying
$\sup_{x\in D} \Vert x - Vx \Vert  < \varepsilon$ and $\Vert V\Vert \leq C$. 
We give several examples of Banach spaces both with and without 
 this approximation property.
Our main results demonstrate that the James-type spaces from a 
general class of quasi-reflexive spaces (which contains 
the classical  James' space $J$) have the
W.A.P,  but that James' tree space $JT$ fails to have the W.A.P.
It is also shown that the dual $J^*$ has the W.A.P.  
It follows that the Banach algebras $W(J)$ and $W(J^*)$, consisting of 
the weakly compact  operators, have
bounded left approximate identities.
Among the other results we obtain a concrete Banach space $Y$ 
so that $Y$ fails to have the W.A.P., but $Y$ has 
this approximation  property without the uniform bound $C$.
} }

\bigskip


\bigskip

\noindent 1. {\bf Introduction.}

\medskip

A Banach space $E$ is said to have the {\it weakly compact approximation
property} (abbreviated W.A.P.)  if there is a constant $C < \infty$
such that  for any weakly compact set $D \subset E$ and $\varepsilon > 0$
there is a weakly compact operator $V: E \to E$ satisfying
$$
\sup_{x\in D} \Vert x - Vx \Vert  < \varepsilon \quad  {\rm and}
\quad\Vert V\Vert \leq C. \eqno{(1.1)}
$$
This (bounded) weakly compact approximation property 
was introduced by Astala and Tylli [AT].
The applications mentioned below were the principal motivation for this in [AT],
but the W.A.P.  is a natural notion worthy of study in its own right.
Clearly any reflexive Banach space has the W.A.P., 
but this property is quite rare for
non-reflexive spaces.  For instance,
if $E$ is a ${\cal L}^1$- or ${\cal L}^{\infty}$-space, then
$E$ has the W.A.P.  if and only if $E$ has the Schur property, 
see [AT,Cor.3]. 
We note that a different notion is obtained by considering
the uniform approximation of the identity operator  on compact sets
by weakly compact operators
(see e.g. Reinov [R], Gr{\o}nb{\ae}k and Willis [GW], 
and Lima, Nygaard and Oja [LNO] for this). 

The weakly compact approximation property defined by (1.1)
has some unexpected applications.  The key fact  [AT,Thm. 1] here is
that the Banach space $F$ has the W.A.P. if and only if 
the measure of weak non-compactness
$$
\omega (S) = \inf \{ \varepsilon > 0: SB_E \subset D + \varepsilon B_F,
D \subset F\ {\rm weakly\ compact} \}
$$
and the weak essential norm
$S \mapsto \Vert S\Vert _w \equiv {\rm dist}(S,W(E,F))$
are uniformly comparable in the space $L(E,F)$ of bounded linear operators
$E \to F$ for all Banach spaces $E$. 
Here  $B_E = \{x \in E: \Vert x\Vert \leq 1\}$ and  $W(E,F)$ stands for
the weakly compact operators $E \to F$. 
The fact that  $c_0$ fails to have the  W.A.P. was then applied  
in [AT,Thm. 4 and Cor. 5]  to show that
$\omega (S)$ is in general  neither uniformly comparable to 
$\omega (S^*)$ nor to $\omega (JS)$ for
arbitrary linear into isometries $J$.
Subsequently, weakly compact approximation properties were
exploited in [T2] to obtain examples of Banach spaces $E$ and $F$, where 
$\Vert S\Vert _w$ is not uniformly comparable to
$\Vert S^*\Vert _w$.  Further applications of the W.A.P. arise from the
fact that the Banach algebra $W(E)$ has a bounded left approximate identity 
whenever $E$ has the W.A.P.

This paper contains several results and examples about  
the weakly compact approximation 
property for Banach spaces.
A first natural question is in which sense "almost reflexive" Banach spaces 
still possess the W.A.P.  A principal aim is to discuss the W.A.P. for 
the class of quasi-reflexive Banach spaces $E$, 
where ${\rm dim}(E^{**}/E) < \infty$.
In sections 2 and 3 we show that the classical
James space $J$ and its dual $J^*$ have the W.A.P.
These results imply that $W(J)$ and
$W(J^*)$ have bounded left approximate identities.
In this direction Loy and Willis [LW] established 
that $W(J)$ has a bounded right approximate identity. 
In section 4 we extend the results of
section 2 by proving that the quasi-reflexive  James-like spaces
constructed by Bellenot, Haydon and Odell in [BHO] have the W.A.P.
These positive results are further highlighted by the 
recent discovery of Argyros and Tolias (see [ArT,Prop. 14.10]) that there 
exist quasi-reflexive hereditarily indecomposable Banach spaces
$E$ that do not have the W.A.P.

In section 5 we present a permanence property for weakly compact approximation 
properties, which implies among other things that certain vector-valued sequence spaces,
including $\ell^1(\ell^p)$ and  $\ell^p(\ell^1)$ for $1 < p < \infty$, have the W.A.P.

Section 6 contains a number of additional examples of  spaces 
failing the W.A.P.
For instance, we show that James' tree space 
$JT$ does not have the W.A.P. We also establish that the W.A.P. differs from the 
corresponding "unbounded" W.A.P. (where the 
uniform bound $\Vert V \Vert \leq C$ is removed from (1.1)). 
Moreover, we obtain a concrete Banach space $Y$ so that the quotient
$Y^{**}/Y$ is isometric to $\ell^2$, but $Y$ does not have the W.A.P.
(This example yields a simpler negative answer to a question from [AT] 
than the quasi-reflexive spaces constructed in [ArT]).
Another natural problem is whether  
$E$ always  has the W.A.P.  if $E$ is  $\ell^1$-saturated and $E$ has a Schauder
basis. Indeed this is false,  since as we show the 
Lorentz sequence spaces $d(w,1)$,
as well as the Azimi-Hagler spaces from [AH], do not have the W.A.P.

\smallskip

The basic terminology and notation related to
Banach spaces will follow  [LT].

\bigskip

\noindent 2. {\bf The James space $J$ has the W.A.P.}

\medskip

The most well-known (and first discovered) quasi-reflexive Banach space $J$ 
was introduced by James [J1]. 
The fact that $J$ has the W.A.P. follows from the more general results
of section 4. However, that argument is more complicated and 
many of the ideas we use there, and in section 3, are well illustrated 
by first presenting them for $J$.

Recall that a real-valued sequence $x = (x_j) \in J$ if 
$\lim_{j\to \infty} x_j = 0$ and the square variation norm
$$
\Vert x \Vert ^2 = \sup \sum_{j=1}^n \vert x_{p_{j+1}} - x_{p_{j}}\vert ^2 < \infty , \eqno{(2.1)}
$$
where the supremum is taken over all indices $1 \leq p_1 < p_2 < \ldots < p_n < p_{n+1}$
and $n \in {\bf N}$.  The monograph  [FG] is a convenient source of results 
(as well as further references)
about $J$. Recall that the coordinate basis $(e_n)$ is a shrinking Schauder
basis for $J$, so that $J^{**}$ can be identified with
the set of scalar sequences $x = (x_j)$ for which
$\sup_n \Vert \sum_{j=1}^n a_je_j\Vert < \infty$. Moreover, 
$J^{**} = \{ x + \lambda {\bf 1}: x \in J, \lambda \in {\bf R}\}$, where ${\bf 1} = (1,1,1,\ldots )$.

The question whether (1.1) is satisfied for any weakly compact subset $D \subset J$ 
can be viewed as
a concrete approximation problem for $J$ that may have independent interest.
The set $\{e_n: n \in {\bf N}\} \cup \{0\}$ of $J$ 
is already a non-trivial test for (1.1), since the sequence
$(e_n)$ is weakly null in $J$.  It turns out that (somewhat surprisingly) the 
desired approximating operators $V \in W(J)$ are 
perturbations of the identity operator by 
certain double averaging functionals over consecutive blocks.

We first state a well known general auxiliary result.
It is convenient to put $[n,m) = \{n,\ldots ,m-1\}$ 
if $m, n \in {\bf N}$ and $m > n$. 

\medskip

\noindent {\sl Lemma 2.1.\   Suppose that $E$ is a Banach space with a normalized Schauder
basis $(e_n)$, and let $D \subset E$ be an arbitrary weakly compact subset.
Then for any $\delta > 0$ and $n \in {\bf N}$ there is $m > n$ such that
for any  $x = \sum _{j=1}^{\infty} a_je_j \in D$ there is an index $j = j(x) \in [n,m)$
satisfying $\vert a_j \vert < \delta$.
}

\medskip

\noindent {\sl Proof.} Suppose to the contrary that there is  $\delta > 0$ and $n \in {\bf N}$
so that for any $m > n$ there is an element $x_m = \sum _{j=1}^{\infty} a_j^{(m)}e_j \in D$
satisfying $\vert a_j^{(m)}\vert \geq \delta$ for all 
$j \in [n,m)$. By the weak compactness of $D$ we may assume that
$x_m \buildrel{w}\over\longrightarrow x = \sum _{j=1}^{\infty} a_je_j \in D$
as $m \to \infty$.
Since $a_j = \lim_{m\to \infty} a_j^{(m)}$ for $j \in {\bf N}$, 
we arrive at  the contradiction that
$\vert a_j \vert \geq \delta$ for each $j \geq n$.
$\square$

\medskip

\noindent {\sl Theorem 2.2. \quad  James' space $J$ has the W.A.P.}

\medskip

\noindent {\sl Proof.}  The argument  will be split into several steps. 
Let $D \subset J$ be a fixed weakly compact subset and $\varepsilon > 0$.
By homogeneity there is no loss of generality to assume that $D \subset B_J$.

\smallskip

\noindent {\it Step 1.} We start by fixing some notation. Given natural numbers 
$1 \leq n_1 < \ldots < n_{k+1}$ we introduce the related averaging functionals  
$A_{[n_{j},n_{j+1})}$ and $A_{(n_{1},\ldots ,n_{k+1})}$ on $J$ by
$$
A_{[n_{j},n_{j+1})}(x) =  {1\over {n_{j+1}-n_{j}}} \sum_{s=n_{j}}^{n_{j+1}-1} a_s, \quad
A_{(n_{1},\ldots ,n_{k+1})}(x) =  {1\over k} \sum_{j=1}^k A_{[n_{j},n_{j+1})}(x),
$$
for  $x = \sum_{i=1}^{\infty} a_ie_i \in J$  and $j =1,\ldots ,k$. Note that 
$\Vert A_{[n_{j},n_{j+1})}\Vert \leq 1$ for each $j$, since
$$
\vert A_{[n_{j},n_{j+1})}(x) \vert = {1\over {n_{j+1}-n_{j}}} 
\vert \langle \sum_{s=n_{j}}^{n_{j+1}-1} e_s^*,x\rangle \vert \leq \Vert x\Vert .
$$
Here $(e_s^*) \subset J^*$ stands for the sequence of biorthogonal functionals 
to $(e_s)$;  clearly
$\Vert e_s^*\Vert = 1$ for  $s \in {\bf N}$.
We next show that certain double averages $ A_{(n_{1},\ldots ,n_{k+1})}$
are uniformly small on the weakly compact set $D$, provided $k$ is large enough.

\smallskip

\noindent  {\sl Claim 1.  Let $\delta > 0$ and $n = n_1 \in {\bf N}$ be arbitrary.
Then there is $k \in {\bf N}$ and natural numbers $n_1 < \ldots < n_{k+1}$ so that
$$
\vert A_{(n_{1},\ldots ,n_{k+1})}(x)\vert < 3\delta \quad {\rm for\ all}\  x \in D.
$$
}

\noindent {\it Proof of Claim 1.} Use Lemma 2.1 repeatedly to choose a sequence 
$n = n_1 < n_2 < \ldots$ in ${\bf N}$,  such that for any $j \in {\bf N}$ and
$x = \sum_{i=1}^{\infty} a_ie_i  \in D$ there is an index $i = i(x)  \in [n_j,n_{j+1})$
satisfying $\vert a_i \vert < \delta$.  Let $k \in {\bf N}$ be given. For each fixed 
$x = \sum_{i=1}^{\infty} a_ie_i \in D$  put
$$
I = \{ j \leq k: {\rm there\ is}\ i \in [n_{j},n_{j+1}) \ {\rm with}\ \vert a_i \vert \geq 2\delta \}
$$
(note that $I$ depends on $x$, $k$ and $\delta$). 
In order to choose $k$ suppose
that $I = \{j_1,\ldots ,j_r\}$ and pick $p_s, q_s \in [n_{j_{s}},n_{j_{s}+1})$
such that $\vert a_{p_{s}}\vert < \delta$ and $\vert a_{q_{s}}\vert \geq 2\delta$
for $s = 1,\ldots ,r$. The square variation norm (2.1) satisfies
$$
\Vert x\Vert \geq (\sum_{s=1}^r \vert a_{q_{s}} - a_{p_{s}}\vert ^2)^{1/2}
\geq \vert I\vert ^{1/2} \delta .
$$
We deduce that the cardinality  $\vert I\vert \leq {1\over {{\delta}^2}}$, 
because $D \subset B_J$ by assumption. 
Note further that 
 $\vert A_{[n_{j},n_{j+1})}(x) \vert 
 < 2\delta$ whenever $j \notin I$, since $\vert a_s\vert < 2\delta$ 
for all $s \in [n_j,n_{j+1})$ in this event. 
By putting these estimates together we get that 
$$
\vert A_{(n_{1},\ldots ,n_{k+1})}(x) \vert \leq {1\over k} (\vert I\vert + (k - \vert I\vert)2\delta )
\leq {1\over {k{\delta}^2}} + 2\delta < 3\delta
$$
once we pick $k > {\delta}^{-3}$. The completes the argument for Claim 1.

\smallskip

\noindent {\it Step 2.}  Fix a decreasing null-sequence $(\varepsilon _j)$ such that
$\sum_{j=1}^{\infty} \varepsilon _j < {\varepsilon}/\sqrt{2}$. 
By successive applications of Claim 1 
we find a sequence of consecutive subdivisions 
$1 = n_{p_{1}} < n_{p_{1}+1} < \ldots < n_{p_{2}}
< n_{p_{2}+1} < \ldots < n_{p_{3}} < n_{p_{3}+1}  < \ldots$ of ${\bf N}$ 
such that 
$$
\vert A_j(x)\vert < \varepsilon _j \quad {\rm for\ all}\ x \in D\ {\rm and}\ j \in {\bf N}, \eqno{(2.2)} 
$$
where we set $A_{j} = A_{(n_{p_{j}},n_{p_{j}+1}, \ldots ,n_{p_{j+1}})}$
for $j \in {\bf N}$.
Let $I_j = [n_{p_{j}},n_{p_{j+1}})$ for $j \in {\bf N}$, which is the "support" 
(with respect to the coordinate basis $(e_s)$) 
of the functional $A_j$ on $J$. 
Put $g_j = \sum_{i\in I_{j}} e_i$, so that $\Vert g_j \Vert \leq \sqrt{2}$ for 
$j \in {\bf N}$.  Define  the linear map $V$ on $J$ by
$$
Vx = x - \sum_{j=1}^{\infty} A_j(x)g_j, \quad x \in J. \eqno{(2.3)}
$$
We verify in three separate steps that $V$ 
provides a uniformly bounded weakly compact approximating operator for
the given weakly compact set $D \subset B_J$ as required by (1.1).

\smallskip

\noindent  {\sl Claim 2. \quad   $\Vert  x - Vx \Vert < \varepsilon$ for all $x \in D$.}

\smallskip

\noindent {\it Proof of Claim 2.}  It follows from (2.3) and the choice of $(\varepsilon _j)$ 
that 
$$
\Vert x - Vx\Vert = \Vert \sum_{j=1}^{\infty} A_j(x)g_j\Vert
\leq \sum_{j=1}^{\infty} \vert A_j(x)\vert \cdot \Vert g_j\Vert < 
{\sqrt{2}} \cdot  \sum_{j=1}^{\infty} \varepsilon _j < \varepsilon .
$$

\smallskip

\noindent  {\sl Claim 3. \quad   $\Vert V \Vert \leq 3$ (independently of the subdivisions).}

\smallskip

\noindent {\it Proof of Claim 3.}  Put $\tilde{V}x =  \sum_{j=1}^{\infty} A_j(x)g_j$ 
for $x \in J$. It suffices to verify that $\Vert \tilde{V}\Vert \leq 2$.

It is convenient to denote $y = \sum_{j=1}^{\infty} y(j)e_j \in J$  
in the argument.  Assume that $x \in J$ is finitely 
supported, and suppose that $q_1 < \ldots < q_{n+1}$ is a sequence of 
coordinates that realizes the square variation norm (2.1) of $\tilde{V}x$. 
We first split
$$
\eqalign{ & \Vert \tilde{V}x\Vert =  
(\sum_{i=1}^n \vert \tilde{V}x(q_{i+1}) - \tilde{V}x(q_i)\vert ^2)^{1/2} \cr
\leq  \ & (\sum_{i\in A} \vert \tilde{V}x(q_{i+1}) - \tilde{V}x(q_i)\vert ^2)^{1/2} +
(\sum_{i\in B} \vert \tilde{V}x(q_{i+1}) - \tilde{V}x(q_i)\vert ^2)^{1/2}, \cr 
}
 \eqno{(2.4)}
$$
where $i \in A$ if both $q_{i}, q_{i+1} \in I_j$  for some $j$, while 
$i \in B$ if $q_{i}$ and $q_{i+1}$ belong to different intervals. 
Note that if  $i \in A$ and $q_{i}, q_{i+1} \in I_j$, then
$\vert \tilde{V}x(q_{i+1}) - \tilde{V}x(q_i)\vert  =  \vert A_j(x) -  A_j(x) \vert  = 0$
 by definition, so that the term 
$(\sum_{i\in A} \vert \tilde{V}x(q_{i+1}) - \tilde{V}x(q_i)\vert ^2)^{1/2}$
actually vanishes.
To estimate the second term
in (2.4) we split
$$
\eqalign{(\sum_{i\in B} \vert \tilde{V}x(q_{i+1}) - \tilde{V}x(q_i)\vert ^2)^{1/2}
\leq \ & (\sum_{i\in B_{1}} \vert \tilde{V}x(q_{i+1}) - \tilde{V}x(q_i)\vert ^2)^{1/2} + \cr 
\ & + (\sum_{i\in B_{2}} \vert \tilde{V}x(q_{i+1}) - \tilde{V}x(q_i)\vert ^2)^{1/2}. \cr
} \eqno{(2.5)}
$$
In (2.5) the set $B_1$ contains every $2k+1$:th term of $B$,  
and $B_2$ the remaining ones.
Consider a single term $\vert \tilde{V}x(q_{i+1}) - \tilde{V}x(q_i)\vert $ 
for some $i \in B_1$, 
and suppose that $q_i \in I_j$ and $q_{i+1} \in I_k$ (where $j < k$). 
In this case 
$$
\vert \tilde{V}x(q_{i}) - \tilde{V}x(q_{i+1})\vert = 
 \vert A_j(x) - A_k(x)\vert .
$$
By definition the double average $A_j(x)$ is  a convex combination of 
$\{x(s): s \in I_j\}$ (and analogously for $A_k(x)$). 
Consequently there are
indices $r_i \in I_j$ and $r_{i+1} \in I_k$ so that
$\vert A_j(x) - A_k(x)\vert \leq \vert x(r_i) - x(r_{i+1})\vert$. 
Since $B_1$ contains every second index from 
$B$, it is easy to check that  the corresponding sequence 
$(r_i)$ is increasing, so we obtain that
$
(\sum_{i\in B_{1}}  \vert A_j(x) - A_k(x)\vert ^2)^{1/2} \leq \Vert x \Vert .
$
By arguing in a similar manner for the sum over the "even" indices $i \in B_2$, 
we get from (2.5) that 
$$
\Vert \tilde{V}x\Vert =  
(\sum_{i\in B} \vert \tilde{V}x(q_{i+1}) - \tilde{V}x(q_i)\vert ^2)^{1/2} \leq 2\Vert x \Vert .
$$
Hence we get by approximating
that $\Vert \tilde{V}x \Vert \leq 2\Vert x\Vert$ for all $x \in J$.
This establishes Claim 3.

\smallskip

\noindent  {\sl Claim 4. \quad   $V \in W(J)$.}

\smallskip

\noindent {\it Proof of Claim 4.} Recall that $V \in W(J)$ if and only if
$V^{**}({\bf 1}) \in J$, where  ${\bf 1} = (1,1,\ldots ) \in J^{**} \setminus J$. 
Put $f_m = \sum_{j=1}^m e_j$ for $m \in {\bf N}$.
Note that $f_{m_{k}} \buildrel{w^*}\over\longrightarrow {\bf 1}$ and
$Vf_{m_{k}} \buildrel{w^*}\over\longrightarrow V^{**}{\bf 1}$ in $J^{**}$
as $k \to \infty$ for each subsequence $(f_{m_{k}})$ of $(f_m)$.
Consider $s_k =  \sum_{j=1}^k g_j$ for $k \in {\bf N}$, 
which determines a subsequence
of $(f_m)$.  Here 
$$
Vg_j = g_j - \sum_{k=1}^{\infty} A_k(g_j)g_k = 0, \quad j \in {\bf N},
$$
since the averages $A_k(g_j) = \delta _{j,k}$ for $j, k \in {\bf N}$.
It follows that $Vs_k = 0$ for $k \in {\bf N}$, so that  $V^{**}{\bf 1} = 0$.
Thus $V \in W(J)$. This completes the proof of Theorem 2.2.
$\square$

\medskip

\noindent {\it Remarks 2.3.}\ (i)\ 
Simpler uniformly bounded approximating operators $V \in W(J)$ 
are available for the particular weakly compact set
$D = \{e_n: n \in {\bf N}\} \cup \{0\} \subset J$. 
Indeed, define $V_k \in W(J)$ for $k \in {\bf N}$ by
$$
V_kx =  x - \sum_{j=1}^{\infty} A_{[kj,k(j+1))}(x)h_j, \quad {\rm where}\ 
h_j = \sum_{s=kj}^{k(j+1)-1} e_s,\  j \in {\bf N}.
$$ 
Then $V_k$ satisfies  (1.1) for $D$ and $\varepsilon > 0$ 
once ${{\sqrt{2}}\over k} < \varepsilon$, since
${1\over {\vert I\vert}} \vert \sum_{s \in I} x(s)\vert \leq {1\over {\vert I\vert}}$
for intervals $I \subset {\bf N}$ and $x = \sum_{s=1}^{\infty} x(s)e_s \in D$.  The 
uniform bound for $\Vert V_k\Vert$ and the weak compactness of $V_k$ 
are easy modifications of Claims 3 and 4 above.

\smallskip

Let $S_k \in L(J)$ be the forward
$k$-shift on $J$ for $k = 0,1,2,\ldots$. 
The reader may also wish to check that  
$V_m = I - {1Ê\over {m+1}} \sum_{k=0}^m S_k \in W(J)$ for $m \in {\bf N}$,
and that $V_m$ satisfies (1.1) for $D$ and $\varepsilon > 0$ 
once $m$ is large enough.

\smallskip

\noindent \ (ii)\ The argument of section 4 yields more complicated
weakly compact approximating operators when applied to $J$
(e.g. the double averaging functionals are replaced 
by less explicit convex combinations, and there is a "shift"-like
perturbation of the identity).

\medskip

The following vector-valued analogue of 
James' space $J$ has been studied in several contexts, 
see e.g.  [PQ]  and  [P].
Let $E$ be a Banach space and $1 < p < \infty$.  
The sequence $x = (x_j) \subset E$
belongs to $J_p(E)$ if $\lim_{j\to\infty} x_j = 0$ and the $p$-variation norm
$$
\Vert x \Vert ^p = 
\sup_{n; p_{1}< \ldots < p_{n+1}} \sum_{j=1}^n \Vert x_{p_{j+1}} - x_{p_{j}}\Vert ^p < \infty .
$$
Here $J_p(E)^{**}/J_p(E) \approx E$ for reflexive spaces $E$,
see e.g. [W,Cor. 2].
Straightforward modifications of the argument in Theorem 2.2 yield that 
$J_p = J_p({\bf R})$ also has the W.A.P.  for $1 < p < \infty$.
More generally, Theorem 5.3 below implies that 
$J_p({\bf R}^n) \approx J_p \oplus \ldots \oplus J_p$ ($n$ summands)
has the W.A.P.  for  $n \in {\bf N}$.
This suggests the following problem (a similar question may obviously be raised for
other James-Lindenstrauss type constructions).

\medskip

\noindent {\it Problem 2.4.}  {\sl Does $J_p(E)$ have the W.A.P. 
whenever $E$ is reflexive and $1 < p < \infty$?
}

\medskip

Recall that a Banach algebra $A$ has a {\it bounded left approximate identity}
(abbreviated B.L.A.I.) if there is a bounded net $(x_{\alpha}) \subset A$
such that 
$$
\lim_{\alpha} \Vert y - x_{\alpha}y\Vert = 0\quad {\rm for}\ y \in A. \eqno{(2.6)}
$$
A  {\it bounded right approximate identity}
(B.R.A.I. for short)  in $A$ is obtained by considering $\Vert y - yx_{\alpha}\Vert$ 
in (2.6). The following observation contains an application of the W.A.P.
to algebras of weakly compact operators.

\medskip

\noindent {\it Proposition 2.5.\  (i)\  If $E$ has the W.A.P., then the Banach algebra 
$W(E)$ has a B.L.A.I. 

\smallskip

\noindent (ii)\  $W(J)$ has a  B.L.A.I.
}

\medskip

\noindent {\it Proof.} (i)\  Let $U \in W(E)$ and $\varepsilon > 0$ be arbitrary.
By applying (1.1) to the weakly compact set $\overline{UB_E}$
we obtain $V \in W(E)$ satisfying $\Vert V \Vert
\leq C$ and
$$
\Vert U - VU\Vert = \sup_{x\in B_{E}} \Vert Ux - VUx\Vert < \varepsilon .
$$
Here $C < \infty$ is a uniform constant. It follows that $W(E)$ 
has a B.L.A.I. by a well-known sufficient condition from [BD,Prop. 11.2].

Part (ii) follows from (i) and Theorem 2.2.
$\square$

\medskip

\noindent {\it Remarks 2.6.} (i)\  Loy and Willis [LW,Cor. 2.4] established 
that $W(J)$ has a B.R.A.I.
using a different matrix-type argument.  Their result suggested the
problem whether $W(J)$ also has a  B.L.A.I.   This question was
not stated in print, but it was known to specialists (R.J. Loy and
G.A. Willis, personal communication). 

Several authors have studied the existence of bounded left or right approximate
identities for algebras of operators, chiefly for closed subalgebras of 
the compact operators on a Banach space. 
For such subalgebras the existence of a B.R.A.I. implies the existence of 
a B.L.A.I. [GW,Cor. 2.7], but no results of this type are known
in the setting of $W(E)$.
We refer to [D,2.9.37 and 2.9.67] for further information and references. 

\smallskip

\noindent (ii)\  The converse of Proposition 2.5.(i) fails 
already for $E = \ell^1$, see [T1, p. 107].

\bigskip

\noindent 3.\quad  {\bf $J^*$ has the W.A.P.}

\medskip

Recall that the dual $J^*$ is also quasi-reflexive of order $1$, but that
$J^*$ is quite different from $J$ as a Banach space. 
For instance, the norm in $J^*$ is not given by any concrete
formula, $J^*$ does not embed into $J$ by [J2,Thm. 3]
and $J$ does not embed into $J^*$ by [A,Thm. 7]
(see also [P] for a different approach).
This provides ample motivation for considering weakly compact approximation 
in $J^*$.

The argument that $J^*$ has the W.A.P.  follows the basic outline of section 2,
but the details are more involved. We recall a few relevant facts about $J^*$.
Put $f_n = \sum_{j=1}^n e_j \in J$ for $n \in {\bf N}$. Then $(f_n)$ is a
boundedly  complete Schauder basis for $J$. 
Put $S^*(x) = \sum_{s \in S} x(s)$
for $x = \sum_{s=1}^{\infty} x(s)f_s \in J$ and any interval $S \subset {\bf N}$.
In the basis $(f_n)$
the square variation norm (2.1) becomes
$$
\Vert x\Vert = \sup_{n;\ S_{1}<\ldots <S_{n}} (\sum_{j=1}^n S_j^*(x)^2)^{1/2},\quad x = \sum_{s=1}^{\infty} x(s)f_s \in J,
\eqno{(3.1)}
$$
where $S_1,\ldots ,S_n$ are intervals of ${\bf N}$ satisfying $\max S_i < \min S_{i+1}$ 
for $i = 1,\ldots ,n-1$ (the interval $S_n$ may be unbounded). We denote this by
$S_1 < S_2 < \ldots < S_n$.
Thus $\Vert S^*\Vert = 1$ whenever $S \subset {\bf N}$ is an interval.  
We will require the fact that
$$
\Vert \sum_{k=1}^n c_kS_k^*\Vert \leq (\sum_{k=1}^n \vert c_k\vert ^2)^{1/2} \eqno{(3.2)}
$$
whenever $S_1 < S_2 < \ldots < S_n$ are intervals of ${\bf N}$ and $c_1,\ldots ,c_n$
are scalars. Indeed,  
$$
\vert \langle \sum_{k=1}^n c_kS_k^*,x\rangle \vert = 
\vert \sum_{k=1}^n c_k S_k^*(x) \vert 
\leq (\sum_{k=1}^n \vert c_k\vert ^2)^{1/2} (\sum_{k=1}^n S_k^*(x)^2)^{1/2}
\leq  (\sum_{k=1}^n \vert c_k\vert ^2)^{1/2} \Vert x\Vert
$$
for $x = \sum_{s=1}^{\infty} x(s)f_s \in J$.

\smallskip

The sequence $(f_n^*)$ of biorthogonal  functionals  to $(f_n)$
forms a $w^*$-basis for $J^*$, that is, 
for any $x^* \in J^*$ there is a unique scalar sequence $(a_j) = (x^*(e_j))$
so that
$x^* = (w^{*}) \sum_{j=1}^{\infty} a_jf_j^*$ as a $w^*$-convergent sum in $J^*$.
It is known that the limit $\lim_{j\to \infty} a_j$ exists for this $w^*$-representation of $x^*$. 
We recall the argument, since we will actually need the  more precise quantitative version 
given below in part (ii).

\medskip

\noindent {\sl Lemma 3.1.  Let $x^* = (w^{*}) \sum_{j=1}^{\infty} a_jf_j^*$ 
be the unique $w^*$-convergent representation of $x^* \in B_{J^*}$.

\smallskip

\noindent\ (i) Then $\lim_{j\to \infty} a_j$ exists.

\smallskip

\noindent\ (ii) There is a uniform constant $C < \infty$
with the following property:
Let $\varepsilon > 0$ and  suppose that for some $k \in {\bf N}$ 
there are  indices $p_1 < q_1 < p_2 < q_2 < \ldots < p_k < q_k$
with $\vert a_{p_{i}} - a_{q_{i}}\vert > \varepsilon$ for $i = 1,\ldots ,k$.
Then 
$
k \leq {C\over {{\varepsilon}^2}}
$.
}

\medskip

\noindent {\it Proof.}  Suppose that there
are $\varepsilon > 0$ and indices $p_1 < q_1 < p_2 < q_2 <  \ldots
< p_k < q_k$
so that $\vert a_{p_{j}} - a_{q_{j}}\vert > \varepsilon$ for all $j \leq k$.
By a result of Casazza, Lin and Lohman (see  [CLL,thm. 16]
or   [FG,Thms. 2.d.1 and 2.c.9])
the sequence $(f_{q_{j}} - f_{p_{j}})_{j=1}^k$ is 
equivalent to the unit vector basis of $\ell^2_k$   with uniform isomorphism 
constants independent of  $(p_j)$,  $(q_j)$ and $k \in {\bf N}$.  
Fix $k$ and consider 
$$
x_k = {1\over {\sqrt{k}}} \sum_{j=1}^k \theta _j(f_{q_{j}} - f_{p_{j}}) \in J,
$$
where the signs $\theta _1,\ldots ,\theta _k$ are chosen so that
$\theta _j (a_{q_{j}} - a_{p_{j}}) = \vert a_{p_{j}} - a_{q_{j}}\vert$ for
$j = 1,\ldots ,k$.
Thus $\Vert x_k \Vert \leq C$, where $C$ is a uniform constant.
We get that
$$
C \geq x^*(x_k) = 
{1\over {\sqrt{k}}} \sum_{j=1}^k \theta _j \langle x^*,f_{q_{j}} - f_{p_{j}}\rangle
= {1\over {\sqrt{k}}} \sum_{j=1}^k \theta _j(a_{q_{j}} - a_{p_{j}}) \geq 
\varepsilon {\sqrt{k}},
$$
which proves both (i) and (ii).
$\square$

\medskip

Suppose that $x^* = (w^{*})  \sum_{j=1}^{\infty} a_jf_j^* \in J^*$ and let $a = \lim_{j} a_j$.
Put  $S_{\infty} = {\bf N}$, so that
$S_{\infty}^*(\sum_{j=1}^{\infty} c_jf_j) = \sum_j c_j$ for 
$\sum_{j=1}^{\infty} c_jf_j \in J$.
We can write
$$
x^* = \sum_{j=1}^{\infty} b_jf_j^* + a S_{\infty}^*, \eqno{(3.3)}
$$
where $b_j \equiv  a_j - a \to 0$ as $j \to \infty$. 
In (3.3) the sum $\sum_{j=1}^{\infty} b_jf_j^*$  is norm-convergent
in $J^*$. 
We will need the following variant of Lemma 2.1 for weakly compact subsets 
of $J^*$.

\medskip

\noindent {\sl Lemma 3.2.  Let $D \subset J^*$ be a
weakly compact set.

\smallskip

\noindent (i) Suppose that $(x_n^*) \subset D$ is a sequence, where
$x_n^* = \sum_{j=1}^{\infty} a^{(n)}_jf_j^* + a_{\infty}^{(n)} S_{\infty}^*$
is written as in (3.3) for $n \in {\bf N}$. 
Then there is a subsequence of $(x_n^*)$, still
denoted by $(x_n^*)$, so that
$$
x_n^*   \buildrel{w}\over\longrightarrow x^* = 
\sum_{j=1}^{\infty} a_jf_j^* + a_{\infty} S_{\infty}^*\quad  {\rm as}\  n \to \infty ,
$$
where $x^*$ is represented as in (3.3) and
$\lim_{n\to\infty} a_j^{(n)} = a_j$ for all $j \in {\bf N} \cup \{\infty\}$. 

\smallskip

\noindent (ii)  For all $n \in {\bf N}$ and $\delta > 0$ there is $m > n$
so that for all $x^* = \sum_{j=1}^{\infty} a_jf_j^* + a_{\infty} S_{\infty}^* \in D$
in the representation (3.3),  there is $j = j(x^*) \in [n,m)$ satisfying 
$\vert a_j \vert < \delta$.
}

\medskip

\noindent {\it Proof.} (i) The weak compactness of $D$ gives
a subsequence of $(x_n^*)$, still denoted by $(x_n^*)$,
so that $x_n^* \buildrel{w}\over\longrightarrow x^*$  as $n \to \infty$. 
Write $x^* =  \sum_{j=1}^{\infty} a_jf_j^* + a_{\infty} S_{\infty}^*$ 
as in  (3.3).
Let $x^{**} \in J^{**}$ satisfy $x^{**}(S_{\infty}^*) = 1$ and
$x^{**}(f_j^*) = 0$ for all $j \in {\bf N}$. Hence 
$a_{\infty}^{(n)} = x^{**}(x_n^*) \to x^{**}(x^*) = a_{\infty}$ as $n \to \infty$.
Moreover,  
$a_j^{(n)} + a_{\infty}^{(n)} = x_n^*(f_j) \to x^*(f_j) = a_j + a_{\infty}$ as 
$n \to \infty$ for $j \in {\bf N}$.
It follows that $\lim_{n\to\infty} a_j^{(n)} = a_j$
for each $j \in {\bf N}$.

\smallskip

\noindent (ii) Suppose to the contrary that $n \in {\bf N}$ and $\delta > 0$
are such that for any $m > n$ there is 
$x_m^* =  \sum_{j=1}^{\infty} a^{(m)}_jf_j^* + a_{\infty}^{(m)} S_{\infty}^* \in D$
represented as in (3.3), for which $\vert a_j^{(m)}\vert \geq \delta$,
for all $j \in [n,m)$. By part (i) there is a subsequence $(x_m^*)$ for which
$$
x_m^* \buildrel{w}\over\longrightarrow x^* = 
\sum_{j=1}^{\infty} a_jf_j^* + a_{\infty} S_{\infty}^* \quad {\rm as}\ m \to \infty ,
$$
where $x^*$ is written as in (3.3) and
$\lim_{m\to\infty} a_j^{(m)} = a_j$ for all $j \in {\bf N}$. 
This implies that $\vert a_j\vert \geq \delta$ for all $j \geq n$, which
contradicts the properties of the expansion (3.3).
$\square$

\medskip

We are ready to prove the main result of this section.

\medskip

\noindent {\sl Theorem 3.3.\quad   $J^*$ has the W.A.P.}

\medskip

\noindent {\it Proof.}  Let $D \subset B_{J^*}$ 
be a fixed weakly compact subset and $\varepsilon > 0$. The desired
approximating operator $V \in W(J^*)$ satisfying (1.1) will again be constructed 
in several stages.

\smallskip

\noindent {\it Step 1.} We first fix some notation. Given $n_1 < \ldots < n_k < n_{k+1}$
we define the averaging functionals
 $B_{[n_{j},n_{j+1})}$ and $B_{(n_{1},\ldots ,n_{k+1})}$ on $J^*$
by
$$
B_{[n_{j},n_{j+1})}(x^*) =  {1\over {n_{j+1}-n_{j}}} \sum_{s=n_{j}}^{n_{j+1}-1} a_s, \quad
B_{(n_{1},\ldots ,n_{k+1})}(x^*) =  {1\over k} \sum_{j=1}^k B_{[n_{j},n_{j+1})}(x^*)
$$
for  $x^* = \sum_{i=1}^{\infty} a_if_i^* + a_{\infty}S_{\infty}^* \in J^*$ written as in (3.3) 
and $j =1,\ldots ,k$.  Note that $B_{[n_{j},n_{j+1})} \in J^{**}$
and that $\Vert B_{[n_{j},n_{j+1})}\Vert \leq 2$ for $j \in {\bf N}$.
In fact, it is easy to see that 
$\vert a_s\vert = \vert x^*(f_s)\vert + \vert a_{\infty}\vert \leq 2\Vert x^*\Vert$
for $x^* = \sum_{i=1}^{\infty} a_if_i^* + a_{\infty}S_{\infty}^* \in J^*$ and $s \in {\bf N}$.

Let $n \in {\bf N}$ and $\delta > 0$ be given. By applying Lemma 3.2.(ii)
repeatedly for $\delta /2 > 0$ we find a sequence $n = n_1 < n_2 < \ldots$ of ${\bf N}$,
such that for any
$x^* = \sum_{i=1}^{\infty} a_if_i^* + a_{\infty}S_{\infty}^* \in D$ and $j \in {\bf N}$
there is $i = i(x^*) \in [n_j,n_{j+1})$ satisfying $\vert a_i\vert < \delta /2$. Here 
$n_{j+1} = m(n_j,\delta /2)$ is given by Lemma 3.2.(ii).
We first verify that sufficiently long double averages 
$B_{(n_{1},\ldots ,n_{k+1})}$ are uniformly small  on 
the weakly compact set $D$.

\smallskip

\noindent {\sl Claim 1.  There is $k \in {\bf N}$ such that
$$
\vert B_{(n_{1},\ldots ,n_{k+1})}(x^*) \vert < 2\delta \quad {\rm for\ all}\ x^* \in D.
$$
}

\noindent {\it Proof.}  The idea resembles that of Claim 1 in Theorem
2.2. Let $x^* \in D$ be arbitrary and write
$x^* = \sum_{i=1}^{\infty} a_if_i^* + a_{\infty}S_{\infty}^*$
as in (3.3). Fix $k \in {\bf N}$ and consider
$$
I = \{j \leq k: \vert a_i\vert \geq \delta \ {\rm for\ some}\ i \in [n_j,n_{j+1})\}
$$
(note that $I$ depends on $x^*, k$ and $\delta$).
Write $I = \{j_1,\ldots ,j_r\}$ and pick $p_s, q_s \in [n_{j_{s}},n_{j_{s}+1})$
such that $\vert a_{p_{s}}\vert \geq \delta$ and $\vert a_{q_{s}}\vert < \delta /2$
for $s = 1,\ldots ,r$.
Since  $\vert a_{q_{s}} - a_{p_{s}}\vert > \delta /2$ for 
each $s = 1,\ldots ,r$,  it follows from Lemma 3.1.(ii)
that $r = \vert I\vert \leq 4C/ {\delta }^2$ for some uniform constant $C < \infty$. 
For $j \notin I$ we have
$\vert a_s\vert < \delta$ for all $s \in [n_j,n_{j+1})$, so that
$\vert B_{[n_{j},n_{j+1})}(x^*) \vert < \delta$. We get the estimates
$$
\eqalign{\vert B_{(n_{1},\ldots ,n_{k+1})}(x^*)\vert \leq \ &
{1\over k} (\sum_{j\in I} \vert B_{[n_{j},n_{j+1})}(x^*) \vert + \sum_{j\notin I} \vert B_{[n_{j},n_{j+1})}(x^*) \vert )
\leq {{2\vert I\vert}\over k} + {{(k - \vert I\vert )\delta}\over k} \cr 
\leq  \ & {{8C}\over {k {\delta}^2}} + \delta  <  2\delta
}
$$
for all large enough $k = k(\delta )$.

\smallskip

\noindent {\it Step 2.} Fix a decreasing positive sequence $(\varepsilon _j)$ such that
$\sum_{j=1}^\infty  \varepsilon _j < \varepsilon$. Next apply Claim 1 successively
to get a sequence of finite subdivisions
$1 = n_{r_{1}} < n_{r_{1}+1} < \ldots < n_{r_{2}} < n_{r_{2}+1} < \ldots < n_{r_{3}} < \ldots$
so that 
$$
\vert B_{(n_{r_{j}},\ldots ,n_{r_{j+1}})}(x^*)\vert < \varepsilon _j 
\ {\rm for\ all}\ x^* \in D \ {\rm and}\ j \in {\bf N}. \eqno{(3.4)}
$$
Put $B_j = B_{(n_{r_{j}},\ldots ,n_{r_{j+1}})}$
and $I_j = [n_{r_{j}},n_{r_{j+1}})$ for $j \in {\bf N}$.
Define the linear map $V$ on $J^*$ by
$$
Vx^* = V(\sum_{i=1}^{\infty} a_if_i^* + a_{\infty} S_{\infty}^*) =
x^* - \sum_{j=1}^{\infty} B_j(x^*)I_j^*, \quad x^* \in J^*.
$$
We next verify that $V$ is a weakly compact operator on $J^*$ 
that satisfies (1.1) for $D$ and the given $\varepsilon > 0$. 
The uniform bound for $\Vert V\Vert$
will require additional tools compared to the argument in section 2. 

\smallskip

\noindent {\sl Claim 2. \quad $\Vert x^* - Vx^* \Vert < \varepsilon$ for all $x^* \in D$.}

\smallskip

\noindent {\it Proof of Claim 2.}  Recall that $\Vert I_j^*\Vert = 1$
for $j \in {\bf N}$ by (3.1). Hence it follows from (3.4) that
$$
\Vert x^* - Vx^*\Vert \leq \sum_{j=1}^{\infty} \vert B_j(x^*)\vert \cdot \Vert I_j^*\Vert 
< \sum_{j=1}^{\infty} \varepsilon _j < \varepsilon , \quad x^* \in D.
$$

\smallskip

\noindent {\sl Claim 3. \quad $\Vert V \Vert \leq 7$ (independently of the subdivisions).} 

\smallskip

\noindent {\it Proof of Claim 3.}
Let $U$ be the linear map 
$$
Ux^* = \sum_{j=1}^{\infty} B_j(x^*)I_j^*, \quad 
x^* \in J^*.
$$
We will write $x^* = \sum_{s=1}^{\infty} x^*(s)f_s^* + x^*_{\infty}S_{\infty}^*  \in J^*$
with respect to the Schauder basis  $\{f_j^*: j \in {\bf N}\} \cup \{S_{\infty}^*\}$ of $J^*$. 
Here $x^*_{\infty} = \lim_{s\to \infty} b_s = \lim_{s\to \infty} x^*(f_s)$
in terms of the $w^*$-representation 
$x^* = (w^{*}) \sum_{s=1}^{\infty} b_sf_s^* \in J^*$. 
Thus $\vert x^*_{\infty} \vert \leq \Vert x^*\Vert$, so that
$$
\Vert \sum_{s=1}^{\infty} x^*(s)f_s^*\Vert \leq \Vert x^*\Vert + \vert x^*_{\infty} \vert \cdot
\Vert S_{\infty}^*\Vert \leq 2\Vert x^*\Vert  \eqno{(3.5)}
$$
for $x^* \in J^*$.  We next establish that 
$$
\Vert Ux^*\Vert \leq 3 \eqno{(3.6)}
$$ 
for all finitely supported $x^* \in [f_1^*,\ldots ,f_n^*] \cap B_{J^*}$,  $n \in {\bf N}$.
Since $[f_1^*,\ldots ,f_n^*] \cap B_{J^*}$ is the closed convex hull of its 
extreme points, it is enough to show that (3.6) holds for the extreme points $x^*$
of $[f_1^*,\ldots ,f_n^*] \cap B_{J^*}$.
Towards this it suffices by Proposition 3.4 below to 
restrict attention to  functionals $x^*$ having the special form 
$$
x^* = \sum_{j=1}^r c_jS_j^*, \quad {\rm where}\ 
\sum_{j=1}^r \vert c_j\vert ^2 = 1\ {\rm and}\  S_1 < \ldots < S_r
$$
are intervals of ${\bf N}$ with $\max S_r \leq n$. Here $n \in {\bf N}$ is arbitrary.
We fix some more notation for convenience.  
Put $I_{j,k} = [n_{r_{j}+k},n_{r_{j}+k+1})$ for $k = 0,\ldots ,r_{j+1} - r_j -1$
corresponding to the subdivision 
$n_{r_{j}} < n_{r_{j}+1} < \ldots < n_{r_{j+1}-1} < n_{r_{j+1}}$
of $I_j = [n_{r_{j}},n_{r_{j+1}})$ for $j \in {\bf N}$.
Thus
${\bf N} = \bigcup_{j=1}^{\infty} (\cup_{k=0}^{r_{j+1}-r_{j}-1} I_{j,k})$
is a partition of ${\bf N}$. 

\smallskip

Let $x^* = \sum_{i=1}^r c_iS_i^*$ be as above. We have 
$$
Ux^* = \sum_{j=1}^{\infty} B_j(x^*)I_j^* = 
\sum_{j=1}^{\infty} (\sum_{i=1}^r c_iB_j(S_i^*))I_j^*
$$
by definition. Note first that  $B_j(S_i^*) = 1$ if the interval
$I_j =  [n_{r_{j}},n_{r_{j+1}}) \subset S_i$. Let 
$A$ consist of the indices $i \in \{1,\ldots ,r\}$ such that 
$I_j \subset S_i$ for some $j \in {\bf N}$, and  put
$\overline{S_i} \equiv \bigcup_{j: I_j \subset S_i} I_j$ for $i \in A$.
Note that $\overline{S_i}$ is an interval contained in $S_i$.
Hence the  corresponding coordinates  satisfy
$Ux^*(s) = c_i$ for $s \in \overline{S_i}$ 
and $i \in A$. 

Suppose next that $I_j = [n_{r_{j}},n_{r_{j+1}}) $ is not contained in
any of the intervals $S_1,\ldots ,S_r$.
For these coordinates  $Ux^*$ looks like
$$
(\sum_{i=1}^r c_iB_j(S_i^*))I_j^* \equiv d_jI_j^*. \eqno{(3.7)}
$$
Here $B_j(S_i^*) = 
{1\over {n_{r_{j+1}}-n_{r_{j}}}} 
\sum_{k=0}^{r_{j+1}-r_{j}-1} {{\vert I_{j,k} \cap S_i \vert}\over {\vert I_{j,k}\vert}} 
\in [0,1)$. Note that  $d_j = \sum_{i=1}^r B_j(S_i^*)c_i$ is an 
"absolutely convex" combination of $c_1,\ldots , c_r$, since
$\sum_{i=1}^r B_j(S_i^*)  \leq 1$.   Put 
$$
E = \{j \in {\bf N}: I_j \nsubseteq S_i\ {\rm for\ any}\ i = 1,\ldots ,r\}.
$$ 
Consequently we may coordinatewise split 
$$
Ux^* = \sum_{i\in A} c_i(\overline{S_i})^* + \sum_{j\in E_{1}} d_jI_j^*
+ \sum_{j\in E_{2}} d_jI_j^* \equiv  \Sigma _1 + \Sigma _2 + \Sigma _3,
$$
where the above sums  are actually finite. Here $E_1$ and $E_2$
contain every second index of $E$, respectively.
It is immediate from (3.2) that   
$\Vert \Sigma _1 \Vert \leq (\sum_{i\in A} \vert c_i\vert ^2)^{1/2} \leq 1$. 

Suppose that $j \in E_1$. By (3.7) one has
$\vert d_j\vert \leq \vert c_{i(j)}\vert  
\equiv \max \{ \vert c_i\vert :  I_{j}\cap S_i \neq  \emptyset \}$ 
for a suitable $i(j) \in \{1,\ldots ,r\}$ with
$I_{j}\cap S_{i(j)} \neq  \emptyset$.
Since $E_1$ contains every second index of $E$
it follows that $i(j_1) \neq i(j_2)$ once $j_1, j_2 \in E_1$
and $j_1 \neq j_2$. 
Thus one gets from (3.2) that
$$
\Vert \Sigma _2\Vert \leq ( \sum_{j\in E_1} \vert d_j\vert ^2)^{1/2}
\leq (\sum_{i=1}^r \vert c_i\vert ^2)^{1/2} = 1.
$$
In a similar manner one checks that  $\Vert \Sigma _3\Vert \leq 1$. 
By putting these estimates together
we get  $\Vert Ux^* \Vert \leq 3$ for these particular functionals
$x^* = \sum_{j=1}^r c_jS_j^* \in [f_1^*,\ldots ,f_n^*] \cap B_{J^*}$, 
which yields (3.6) (modulo Proposition 3.4 below). 

Finally, from (3.5), (3.6) and $U(S_{\infty}^*) = 0$ we obtain by approximation that
$$
\Vert Uy^* \Vert = \Vert U(\sum_{s=1}^{\infty} y^*(s)f_s^*)\Vert
\leq 3\Vert \sum_{s=1}^{\infty} y^*(s)f_s^*\Vert \leq 
6\Vert y^*\Vert 
$$
for $y^* = \sum_{s=1}^{\infty} y^*(s)f_s^* + y_{\infty }^*S_{\infty}^* \in J^*$.
We deduce that $\Vert V \Vert = \Vert I - U \Vert  \leq 7$.

\smallskip

\noindent {\sl Claim 4. \quad $V \in W(J^*)$.}

\smallskip

\noindent {\it Proof of Claim 4.} Let $(x_n^*) \subset B_{J^*}$ be an arbitrary 
sequence. We are required to find a subsequence, still denoted by
$(x_n^*)$, so that $(Vx_n^*)$ is weakly convergent.  
Write $x_n^* = \sum_{j=1}^{\infty} a_j^{(n)}f_j^* + a_{\infty}^{(n)}S_{\infty}^*$
as in (3.3) for $n \in {\bf N}$.
By the $w^*$-sequential compactness and a diagonalization we may pass 
to a subsequence and assume without loss of generality
that
$\lim_{n\to\infty}  a_j^{(n)} = a_j$ for  $j \in {\bf N} \cup \{\infty \}$,
$\lim_{j\to \infty} a_j = a$, and
$$
x_n^*  \buildrel{w^*}\over\longrightarrow x^* \equiv (w^*)
\sum_{j=1}^{\infty} a_jf_j^* + a_{\infty}S_{\infty}^* =
\sum_{j=1}^{\infty} (a_j - a)f_j^* + (a_{\infty} + a)S_{\infty}^* 
\in J^* 
$$
as $n \to \infty$. 
The latter representation is the norm-convergent one as in (3.3). Put
$$
y_n^* = x_n^* - a \sum_{j=1}^n I_j^* = 
\sum_{s=1}^{\max I_n}  (a_s^{(n)} - a)f_s^* + 
\sum_{s= \max I_{n}+1}^{\infty} a_s^{(n)}f_s^*
+ a_{\infty}^{(n)}S_{\infty}^*
$$
as a norm-convergent sum for $n \in {\bf N}$.
Let $S_{\infty}^{**} = w^* -\lim_{n\to \infty} f_n$ in $J^{**}$, so that
$S_{\infty}^{**}(S_{\infty}^*) = 1$ and $S_{\infty}^{**}(f_j^*) = 0$
for $j \in {\bf N}$. Put $g_j = f_j - S_{\infty}^{**} $ for $j \in {\bf N}$.
Then $\{g_j: j \in {\bf N}\} \cup \{S_{\infty}^{**}\}$ is a Schauder 
basis for $J^{**}$, for which 
$$
\lim_{n\to\infty} g_j(y_n^*) = \lim_{n\to\infty} (a_{j}^{(n)} - a) = a_j - a\quad 
(j \in {\bf N}),\
\lim_{n\to\infty} S_{\infty}^{**}(y_n^*) = 
\lim_{n\to\infty} a_{\infty}^{(n)} = a_{\infty}.
$$
This yields that 
$
y_n^*  \buildrel{w}\over\longrightarrow y^* \equiv \sum_{j=1}^{\infty} (a_j - a)f_j^*
+ a_{\infty}S_{\infty}^*
$
(norm-convergent sum) in $J^*$ as $n \to \infty$.

Note further that $V(I_j^*) = I_j^* - \sum_{k=1}^{\infty} B_k(I_j^*)I_k^*
= 0$ for $j \in {\bf N}$, since $B_k(I_j^*) = \delta _{k,j}$ 
for $j, k \in {\bf N}$ by definition. 
Thus $Vx_n^* = Vy_n^*  \buildrel{w}\over\longrightarrow Vy^*$  in $J^*$ as
$n \to \infty$.
The proof of Theorem 3.3 will be complete once we have verified the
following auxiliary fact that was used in the proof of Claim 3.
$\square$

\medskip

\noindent {\sl Proposition 3.4.\   Put $J_n^* = [f_1^*,\ldots ,f_n^*]
\subset J^*$ for $n \in {\bf N}$.
Then the extreme points of $B_{J_{n}^*}$ are contained in the set of elements
of the form 
$$
\sum_{j=1}^k c_jS_j^*,\quad {\rm where}\
\sum_{j=1}^k \vert c_j\vert ^2 = 1\ {\rm and}\ S_1 < \ldots < S_k
$$
are intervals of ${\bf N}$ with $\max S_k \leq n$.
}

\medskip

\noindent {\it Proof.}  It follows from (3.1) and (3.2) that 
$$
D = \{\sum_{j=1}^k c_jS_j^*: S_1< \ldots < S_k\ {\rm intervals,}\ 
\max S_k \leq n\ {\rm and}\ \sum_{j=1}^k \vert c_j\vert ^2 = 1\} \subset B_{J_{n}^*}
$$
is a symmetric $1$-norming set for $J_n = [f_1,\ldots ,f_n] \subset J$. 
Indeed,  it is immediate from (3.2) that
$\vert (\sum_{j=1}^kc_jS_j^*)(x)\vert \leq \Vert x\Vert$ for
$\sum_{j=1}^k c_jS_j^* \in D$ and $x \in J_n$.
Conversely, for any non-zero 
$x \in J_n$ there are intervals $S_1 < \ldots < S_k$ with $\max S_k \leq n$,
so that $\Vert x \Vert ^2 = \sum_{j=1}^k S_j^*(x)^2$. By choosing 
$c_j = {\Vert x\Vert }^{-1} S_j^*(x)$ for $j = 1,\ldots ,k$ we get that
$$
\Vert x \Vert = \sum_{j=1}^k c_jS_j^*(x) = (\sum_{j=1}^kc_jS_j^*)(x),
$$
where $\sum_{j=1}^k \vert c_j\vert ^2 = 1$. 

Observe next that $\overline{co}(D) = B_{J_{n}^*}$. In fact, if 
$\overline{co}(D) \varsubsetneq B_{J_{n}^*}$, then the 
Hahn-Banach theorem would give $x_0^*
 \in B_{J_{n}^*}$ and $x \in B_{J_n}$ so that $x_0^*(x) = 1$ and
$x^*(x) \leq \alpha < 1 = \Vert x\Vert$ for all $x^* \in \overline{co}(D)$.
This contradicts the $1$-norming property of $D$.
Finally, since $B_{J_{n}^*} = \overline{co}(D)$, Milman's "converse" 
to the Krein-Milman theorem
(see e.g. [Ph,Prop. 1.5]) yields that the set of extreme points 
$ext(B_{J_{n}^*}) \subset D$.  $\square$

\medskip

Proposition 2.5.(i) and Theorem 3.3 have the following consequence. 

\medskip

\noindent {\sl Corollary 3.5.\quad   $W(J^*)$ has a
B.L.A.I.}

\medskip

\noindent {\sl Remarks 3.6.} \ (i)  Recall that $J^{**} \approx J$, so that
$J^{***} \approx J^*$. Hence Theorems 2.2 and 3.3 imply that 
the $k$:th dual $J^{(k)}$ has the W.A.P.  for all $k \geq 2$.
Moreover, $J = [f_n^*: n \in {\bf N}]^*$, 
since $(f_n)$ is a monotone boundedly complete basis for $J$. 
It follows that the predual $J_* \equiv 
[f_n^*: n \in {\bf N}]$ of $J$ also has the  W.A.P., 
since $J^* \approx J_* \oplus [S_{\infty}^*]$ and 
the W.A.P. is
inherited by complemented subspaces (cf. Lemma 5.2.(i) below). 

\smallskip

\noindent (ii)\  The averaging functionals $B_j \in J^{**}$ employed 
in the proof of Theorem 3.3 are not $w^*$-continuous on $J^*$, so that
the approximating operators $V \in W(J^*)$ are not adjoints of operators on $J$.
Hence Corollary 3.5 does not (by itself) imply  the earlier result of 
[LW,Cor. 2.4]  that  $W(J)$ has a B.R.A.I.

\bigskip

\noindent 4.\quad  {\bf A family of James-type spaces having the W.A.P.}

\medskip

The class of quasi-reflexive Banach spaces is  extensive, and 
sections 2 and 3 suggested the question whether 
there are quasi-reflexive spaces $E$ that {\it fail} to have the W.A.P.
During the course of this work Argyros and Tolias discovered
that a class of hereditarily indecomposable (H.I.) spaces constructed 
recently (for different purposes) in [ArT]  contains quasi-reflexive spaces
of this kind.
We refer to [ArT] for the description of these
spaces, and to [ArT,Prop. 14.10] for the details of the following example. 

\medskip

\noindent {\sl Example 4.1.\ There 
is a quasi-reflexive H.I. space $E$ that fails to have the W.A.P.
}

\medskip

One reason for such spaces $E$ 
to fail the W.A.P. appears to be that 
they admit "few" weakly compact operators in the sense that
$$
L(E) = \{ \lambda I + V: \lambda \in {\bf C}, V\ {\rm is\ strictly\ singular\
and\ weakly\ compact}\}.
$$
In contrast the quasi-reflexive spaces studied in this paper
have many weakly compact subsets, but they are also 
sufficiently rich in weakly compact operators to suggest that they may have
the W.A.P.
The purpose of this section is to extend the results of section 2 
to a general class of quasi-reflexive spaces considered 
by Bellenot, Haydon and Odell  [BHO]. 
This class also contains $J$, but the reader is expected to 
already be familiar with the argument from section 2. 
The desired approximating operators will also be somewhat 
more involved in the general case.

\smallskip

Let $(h_j)$ be a normalized Schauder basis for a reflexive space $E$.
The Banach space $J(h_j)$ consists of the scalar sequences
$(a_j)$ so that $\lim_{j\to\infty} a_j = 0$ and
$$
\Vert (a_j)\Vert = \sup \{ \Vert \sum_{j=1}^n (a_{p_{j}} - a_{q_{j}})h_{p_{j}} \Vert : 
1 \leq p_{1} < q_{1} < \ldots < p_{n} < q_{n}, n \in {\bf N}\} < \infty .
\eqno{(4.1)}
$$
We obtain $J$ with an equivalent norm to (2.1),
if $(h_j)$ is the standard coordinate basis of $\ell^2$.
The reference [BHO] contains the basic information about this construction,
where it is discussed in terms the boundedly complete basis 
(analogous to (3.1) for $J$).
All these  spaces $J(h_j)$ are quasi-reflexive of order $1$ by [BHO,Thm. 4.1].
In addition, $J(h_j) \approx J(u_j) \approx J(g_j)$ [BHO,Prop. 1.1],
where $(u_j)$ is the {\it unconditionalization} of $(h_j)$ defined by 
$$
\Vert \sum_{j=1}^{\infty} a_ju_j\Vert =
\sup \{ \Vert \sum_{j=1}^{\infty} \theta _ja_jh_j\Vert : (\theta _j) \in \{-1,1\}^{\bf N}\} < \infty,
$$
and $(g_j)$ is the {\it right dominant} version of $(h_j)$ (or of $(u_j)$)
given by 
$$
\Vert \sum_{j=1}^{\infty} a_jg_j\Vert =
\sup \{ \Vert \sum_{j=1}^{\infty} a_{n(i)}h_{m(i)}\Vert : 
1 \leq m(1) \leq n(1) < m(2) \leq n(2) < \ldots \} < \infty .
$$
Hence we may and will assume in the sequel that the original basis $(h_j)$ of $E$
is $1$-unconditional,
and in view of  [BHO,Prop. 1.1.(4)] that
$$
\Vert \sum_{i=1}^{\infty} a_{n(i)}h_{m(i)}\Vert \leq 
2\Vert \sum_{i=1}^{\infty} a_{n(i)}h_{n(i)}\Vert  \eqno{(4.2)}
$$
for all $\sum_{j=1}^{\infty}  a_jh_j \in E$ and all sequences
$1 \leq m(1) \leq n(1) < \ldots < m(i) \leq n(i) < \ldots$.

Let $(e_j)$ stand for the unit vector basis in $J(h_j)$, which is a normalized
monotone Schauder basis for $J(h_j)$. 
Recall that the basic sequence 
$(x_j)$ in $J(h_j)$ is a {\it skipped block} basic sequence of $(e_j)$
if for all $j$ there is $n(j) \in {\bf N}$ so that 
$\max supp(x_j) < n(j) < \min supp(x_{j+1})$.  
Here, as well as in the sequel, the support $supp(x)$ of $x \in J(h_j)$ 
is with respect to the basis $(e_j)$.
Every normalized skipped block basic sequence 
$(x_n)$ of $(e_j)$ is $C$-unconditional 
in $J(h_j)$ with a uniform constant $C < \infty$, see [BHO,Prop. 2.1.(2)].

\smallskip

We will require some additional tools and auxiliary results 
in order to extend the argument of section 2 to the present setting.
Our first result concerns the uniform unconditionality of
skipped block sequences of a special type in $J(h_j)$.
It is convenient to denote the natural projection of $J(h_j)$ 
onto the span $[e_s: m \leq s < n]$ by $P_{[m,n)}$,
that is,  $P_{[m,n)}=  P_{n-1}(I - P_{m-1})$ for $2 \leq m < n$.
Here $(P_n)$ are the basis projections on $J(h_j)$ with respect to
$(e_j)$.
Thus $\Vert P_{[m,n)}\Vert  \leq 2$ for $m \geq 2$.

\medskip

\noindent {\sl Lemma 4.2.\ Let $(x_j)$ be a skipped block basic sequence of $(e_j)$
and let  $S_j \subset {\bf N}$ be the smallest intervals 
satisfying $S_j \supset supp(x_j)$ for $j \in {\bf N}$.
Suppose that $(n_j) \subset {\bf N}$ is an increasing  sequence so that
$n_0 = 1$ and $n_j \in S_{j+1}$ for $j \in {\bf N}$, and  put 
$$
y_j = P_{[n_{j-1},n_{j})}(x_j + x_{j+1}), \quad j \in {\bf N}.
$$
Then there is an absolute
constant $C < \infty$ so that
$$
\Vert \sum_{j\in B} a_jy_j\Vert \leq 
C\Vert \sum_{j=1}^{\infty} a_jy_j\Vert  \eqno{(4.3)}
$$
for all subsets $B \subset  {\bf N}$ and all norm convergent
$\sum_{j=1}^{\infty} a_jy_j \in J(h_i)$. 
}

\medskip

\noindent {\it Proof.}  Let $B \subset {\bf N}$ be a given set. By approximation 
it is enough to establish (4.3) with a uniform constant $C$ for all finitely
supported sums $\sum_j a_jy_j$. Put $z = \sum_{j\in B} a_jy_j$ and $\tilde{z} = 
\sum_j a_jy_j$. Suppose that 
$p_{1} < q_{1} < \ldots < p_{m} < q_{m}$
is a sequence of coordinates norming $z$, so that 
$$
\Vert z \Vert =  \Vert \sum_{j=1}^m  (z(p_{j}) - z(q_{j}))h_{p_{j}} \Vert  = 
\Vert \sum_{j=1}^m  b_jh_{p_{j}} \Vert ,
$$
where we put  $b_j = z(p_{j}) - z(q_{j}) \neq 0$ for  $j = 1,\ldots ,m$.
Here we use the convenient notation $y = \sum_{s=1}^{\infty} y(s)e_s$
for elements $y \in J(h_j)$.

Let $T_r \subset {\bf N}$ be the smallest interval satisfying
$T_r \supset supp(y_r)$ for $r \in {\bf N}$. We put 
$E = \{j \leq m:\   p_j, q_j \in T_r$ for some $r\}$ and 
$F = \{j \leq m:\   p_j \in T_r, q_j \in T_s$ for $r  <  s\}$. 
We initially split
$$
 \Vert \sum_{j=1}^m b_jh_{p_{j}} \Vert \leq \Vert \sum_{j\in E} b_jh_{p_{j}}\Vert +
\Vert \sum_{j\in F_{1}} b_jh_{p_{j}}\Vert +  
\Vert \sum_{j\in F_{2}} b_jh_{p_{j}}\Vert 
\equiv \Sigma _1 + \Sigma _2 + \Sigma _3. \eqno{(4.4)}
$$
Here $F = F_1 \cup F_2$ is the partition of $F$ into every second index.
Observe first that $\tilde{b}_j \equiv \tilde{z}(p_j) - \tilde{z}(q_j) = z(p_j) - z(q_j) = b_j$
whenever $j \in E$, so that $\Sigma _1 \leq \Vert \tilde{z}\Vert$ by definition.

We next verify that $\Sigma _2 \leq 2\Vert \tilde{z}\Vert$. 
Let $j \in F_1$. There are three cases to consider:
$$
(i)\  p_{j} \in T_r,\ q_{j} \in T_s,\quad
(ii)\ p_{j} \in T_s,\  q_{j} \in T_t,\quad (iii)\  p_{j} \in T_r,\ q_{j} \in T_t, \eqno{(4.5)}
$$
where $r < s < t$, $r, t \in B$ and $s \notin B$. 
In case (i) from (4.5) one has $b_{j} = z(p_{j}) - z(q_{j}) = 
z(p_{j})$.  
Since $(x_s)$ is a skipped block sequence, 
and $y_s = P_{[n_{s-1},n_{s})}(x_s + x_{s+1})$,
there is some index $r_{j} \in T_s$ for which $y_s(r_{j}) = 0$.
Hence we may move $q_{j}$ to $r_{j}$, so that the difference 
$$
b_{j} = z(p_{j}) = \tilde{z}(p_{j}) -  \tilde{z}(r_{j}) \equiv \tilde{b}_{j}
$$
can still be used towards computing $\Vert \tilde{z}\Vert$ as in 
definition (4.1). 
For case (ii) from (4.5) observe first that 
$b_{j} = z(p_{j}) - z(q_{j}) = - z(q_{j})$. 
There are  two possibilities to consider.
If $\vert \tilde{z}(p_{j}) - \tilde{z}(q_{j})\vert  > 
{{\vert b_{j}\vert}\over 2}$, then we keep 
the coordinates $p_{j} < q_{j}$
towards computing $\Vert \tilde{z}\Vert$ as in (4.1).
In the opposite case, where $\vert \tilde{z}(p_{j}) - \tilde{z}(q_{j})\vert  \leq 
{{\vert b_{j}\vert }\over 2}$, we move the coordinate $q_{j}$
to some $r_{j} \in T_t$ satisfying $y_t(r_{j}) = 0$. This implies that
$$
\vert \tilde{z}(p_{j}) - \tilde{z}(r_{j})\vert 
= \vert \tilde{z}(p_{j})\vert \geq \vert z(q_{j})\vert 
- \vert z(q_{j}) - \tilde{z}(q_{j})\vert = 
\vert b_j\vert - \vert \tilde{z}(q_{j}) - \tilde{z}(q_{j})\vert
\geq {{\vert b_{j}\vert}\over 2}.
$$
Finally, in case (iii) we have 
$z(p_j) - z(q_j) = \tilde{z}(p_j) - \tilde{z}(q_j)$ (because $r, t \in B$),
and we retain the pair $p_j < q_j$.

We get in all cases from (4.5) that
$\vert b_{j}\vert \leq 2\vert \tilde{b}_{j}\vert$,
where $\tilde{b}_{j} \equiv \tilde{z}(p_{j}) - \tilde{z}(\tilde{q}_{j})$
and $\tilde{q}_{j}$ stands for
either $q_{j}$ or $r_{j}$, depending on the indicated choices.
The $1$-unconditionality of the basis $(h_i)$ in $E$ yields then that
$$
\Sigma _2 = \Vert \sum_{j\in F_{1}} b_jh_{p_{j}}\Vert \leq 
2 \Vert \sum_{j\in F_{1}} \tilde{b}_jh_{p_{j}}\Vert \leq 2\Vert \tilde{z}\Vert ,
$$
since the sequence of pairs  $p_{j} < \tilde{q}_{j}$ with $j \in F_1$
are admissible coordinates towards computing $\Vert \tilde{z}\Vert$ as in (4.1). 
Indeed, the fact that $F_1$
contains every second index of $F$ ensures that the order
is preserved in the new sequence if (some of) the coordinates $q_{j}$ are moved. 

The estimate $\Sigma _3 \leq 2\Vert \tilde{z}\Vert$ is  similar.
This completes the proof of (4.3).
$\square$

\medskip

The following combinatorial lemma due to Ptak [Pt]  (see also 
[BHO] for the present formulation) will be a crucial tool
towards building certain convex combinations of averaging functionals, 
which will replace the double averages used in sections 2 and 3.
It is convenient to put $(\alpha _1,\ldots ,\alpha _k) \in (S_{\ell^k_{1}})_{+}$ 
provided $\sum_{r=1}^k \alpha _r = 1$ and  $\alpha _r \geq 0$ for $r = 1,\ldots ,k$.

\medskip

\noindent {\sl Lemma 4.3.\ } [BHO,Lemma 3.1] {\sl Let $0 < \delta < 1$ be fixed.
Suppose that ${\cal F}$ is a collection of non-empty finite subsets
of ${\bf N}$ satisfying the following properties:

\smallskip

\noindent \ (i)\ $B \in {\cal F}$ whenever 
$\emptyset \neq B \subset A$ and $A \in {\cal F}$.

\smallskip

\noindent \ (ii)\  For every $k \in {\bf N}$ and every convex combination
$(\alpha _1,\ldots ,\alpha _k) \in (S_{\ell^k_{1}})_{+}$, 
there is $A \in {\cal F}$ so that
$\sum_{j\in A} \alpha _j \geq \delta$.

\smallskip

\noindent Then there is an infinite subset $M \subset {\bf N}$ so that
$A \in {\cal F}$ for all non-empty finite sets $A \subset M$.
}

\medskip

We are now ready to prove the main result of this section, which 
establishes the weakly compact approximation property
for the James-like spaces $J(h_j)$.
In the proof we will denote $x \in J(h_j)$ by $x = \sum_{s=1}^{\infty} x(s)e_s$.

\medskip

\noindent {\sl Theorem  4.4.\  Let $(h_j)$ be a normalized Schauder  basis
for a reflexive Banach space $E$. 
Then the quasi-reflexive space $J(h_j)$ has the W.A.P.
}

\medskip

\noindent {\it Proof.} Recall that in view of [BHO,Prop. 1.1]  we 
may assume that the Schauder basis $(h_j)$ for $E$ is $1$-unconditional
and satisfies the right dominance property (4.2).
Let $D \subset B_{J(h_{j})}$ be a weakly compact set
and let $\varepsilon > 0$.
We again split the argument into 
distinct steps.

\smallskip

\noindent {\it Step 1.} Let $0 < \delta < 1$ and $n \in {\bf N}$ be given.  
By successive applications of Lemma 2.1 we fix
a sequence $n = n_1 < n_2 < \ldots$ in ${\bf N}$, so that 
for every $x = \sum_{s=1}^{\infty} x(s)e_s \in D$
and $j \in {\bf N}$ there is some index $s_j \in [n_{j},n_{j+1})$ for which
$\vert x(s_j)\vert < {{\delta}\over {2^{j+3}}}$. 

\smallskip

For technical purposes we need to improve 
this fact by a slight perturbation. 

\medskip

\noindent {\sl Claim 1.\   For every $x \in D$ there is a perturbation 
$x \approx x_0 + \sum_{j=1}^{\infty} a_jy_j$
satisfying the following properties:

\smallskip

\item{(i)}$x_0 = P_{[1,n_{1})}(x)$,
$y_j \in [e_s:  n_j \leq s < n_{j+1}]$, $\Vert y_j\Vert = 1$
and $a_j \geq 0$ for $j \in {\bf N}$,

\smallskip

\item{(ii)} $\Vert x - (x_0 + \sum_{j=1}^{\infty} a_jy_j)\Vert < {{\delta}\over {4}}$
and $\Vert \sum_{j=1}^{\infty} a_jy_j\Vert < 2 + {{\delta}\over 4} < 3$, 

\smallskip

\item{(iii)}  for every $j \in {\bf N}$ one has $y_j(s_j) = 0$ 
for some $s_j \in [n_j,n_{j+1})$.
}

\medskip

\noindent {\it Proof of Claim 1.} Let $x = \sum_{s=1}^{\infty} x(s)e_s \in D$.
We put $x_0 = P_{[1,n_{1})}(x)$ and 
$$
v_j = \sum_{s\in [n_{j},n_{j+1}); s \neq s_{j}} x(s)e_s \in
[e_s: n_j \leq s <n_{j+1}, s \neq s_{j}], \quad j \in {\bf N}.
$$
Consider $y = x_0 + \sum_{j=1}^{\infty} v_j = x_0 + \sum_{j=1}^{\infty} a_jy_j$,
where $a_j = \Vert v_j\Vert$ and $y_j = {{v_j}\over {\Vert v_{j}\Vert}}$
(with obvious modifications if $v_j = 0$).
Let $z_j = P_{[n_{j},n_{j+1})}(x)$, so 
that $\Vert z_j - v_j\Vert = \vert x(s_j)\vert < 
 {{\delta}\over {2^{j+3}}}$ for $j \in {\bf N}$. Hence
$
\Vert x - y\Vert \leq  \sum_{j=1}^{\infty} \Vert z_j - v_j\Vert 
< {{\delta}\over 4}$.
The other conditions from Claim 1 are satisfied by construction.

\smallskip

We note that in Claim 1 we also get that
$$
a_j \leq \Vert P_{[n_{j},n_{j+1})}(x)\Vert  + \Vert z_j - v_j\Vert <
2 +  {{\delta}\over {2^{j+3}}} < 3,\quad j \in {\bf N}. \eqno{(4.6)}
$$

\smallskip

\noindent {\it Step 2.}  Define the averaging functionals $A_j$ on $J(h_i)$ by
$$
A_j(x) = {1\over {n_{j+1} - n_{j}}} \sum_{s=n_{j}}^{n_{j+1}-1} x(s) \quad {\rm for}\ 
x = \sum_{s=1}^{\infty} x(s)e_s \in J(h_j)
$$
and $j \in {\bf N}$. 
Thus $\Vert A_j \Vert = 1$ for $j \in {\bf N}$, since 
$\vert x(s)\vert \leq \Vert x\Vert$ for $x \in J(h_j)$ and $s\in {\bf N}$.
For any given  $k \in {\bf N}$ and $(\alpha _1,\ldots ,\alpha _k) \in (S_{\ell^k_{1}})_{+}$
we introduce  the convex combination  
$$
A(x) = \sum_{j=1}^k \alpha _jA_j(x),\quad x \in J(h_j).  \eqno{(4.7)}
$$
By definition  $A(x)$ is a convex combination 
of the coordinates $\{x(s): s =n_1,\ldots ,n_{k+1}-1\}$
of $x = \sum_{s=1}^{\infty} x(s)e_s \in J(h_j)$, and $\Vert A \Vert \leq 1$.
Note that $A$ depends on 
$(n_j)$, $k$ and  $(\alpha _1,\ldots ,\alpha _k)$, 
but  our notation does not make this explicit for simplicity.

The double averages involved in the weakly compact approximating operators
from sections 2 and 3  relied implicitly on the square variation norm for $J$. 
Our next aim is to show (see Claim 2 below) that 
we may choose $k \in {\bf N}$ and a convex combination
$(\alpha _1,\ldots ,\alpha _k) \in (S_{\ell^k_{1}})_{+}$,
so that the corresponding $A$ from (4.7) satisfies
$\vert A(x)\vert < \delta$ for every $x \in D$.
For this end  Lemma 4.3 will be crucial. We
formulate the main technical step here as a separate lemma.

\medskip

\noindent {\sl Lemma 4.5.\  There is an integer $k \in {\bf N}$ 
and a convex combination
$(\alpha _1,\ldots ,\alpha _k) \in (S_{\ell^k_{1}})_{+}$ so that for
all $x \in D$ and perturbations $x \approx x_0 + \sum_{i=1}^{\infty} a_iy_i$ as in Claim 1, 
one has
$$
\sum_{i=1}^k  \alpha _ia_i < {{\delta}\over {4}}.  \eqno{(4.8)}
$$
}

\noindent {\it Proof of Lemma 4.5.} 
Assume to the contrary that for every $k \in {\bf N}$ and every 
convex combination $(\alpha _1,\ldots ,\alpha _k) \in (S_{\ell^k_{1}})_{+}$
we may find an element  $x \in D$ and a 
perturbation $x \approx x_0 + \sum_{j=1}^{\infty} a_jy_j$
satisfying conditions (i) - (iii) of Claim 1 and (4.6), 
so that 
$$
\sum_{j=1}^k  \alpha _ja_j \geq  {{\delta}\over {4}}.  \eqno{(4.9)}
$$

We wish to apply the combinatorial Lemma 4.3 to this setting.
Let ${\cal F}$ be the collection of finite sets 
$$
\eqalign{A = \{i:   \ a_i \geq \ & {{\delta}\over 8},  {\rm where}\ 
x_0 + \sum_{j=1}^{\infty} a_jy_j\ {\rm satisfies}\ (i), (ii), (iii)
\ {\rm of\ Claim\ 1}, (4.6), \cr  
\ & {\rm and} \ (4.9)\ {\rm holds\ for\ some}\ k\ {\rm and}\ 
(\alpha _1,\ldots ,\alpha _k) \in (S_{\ell^k_{1}})_{+}\}. \cr
}
$$
The family ${\cal F}$ satisfies the conditions of Lemma 4.3. Indeed, 
$$
3\cdot \sum_{i\in A} \alpha _i \geq \sum_{i\in A} \alpha _ia_i \geq  
{{\delta}\over 4} - \sum_{i \in A^c} \alpha _ia_i  > {{\delta}\over 8},
$$
since $\sum_{i \in A^c} \alpha _ia_i < {{\delta}\over 8}$.
Here $A^c = \{1,\ldots ,k\} \setminus A$.

Lemma 4.3 yields an infinite set 
$M = \{m_i: i \in {\bf N}\} \subset {\bf N}$
for which $A \in {\cal F}$ for all finite subsets $A \subset M$.
By applying this fact successively to 
$A_n = \{m_1,\ldots ,m_n\}$ for $n \in {\bf N}$,
we obtain  a sequence of elements 
$z_n = x_0^{(n)} + \sum_{j=1}^{\infty} a_j^{(n)}y_j^{(n)} \in J(h_j)$ 
satisfying conditions (i) - (iii) of Claim 1 and (4.6), and where further
$$
a^{(n)}_{m_{i}} \geq {{\delta}\over 8} \quad {\rm for}\  i = 1,\ldots ,n. \eqno{(4.10)}
$$
By a compactness argument we may pass to a subsequence of $(z_n)$, so that
$a_i^{(n)} \to a_i$, $y_i^{(n)} \to v_i$ and $x_0^{(n)} \to v$
 in norm for every  $i \in {\bf N}$ as 
$n \to \infty$. Here the supports satisfy $supp(v) \subset [1,n_1)$, 
$supp(v_i) \subset [n_i,n_{i+1})$,
and by condition (iii) we may further ensure that $v_i(s_i) = 0$ for
some $s_i \in [n_i,n_{i+1})$ for $i \in {\bf N}$.

Lemma 4.2  implies that the sequence 
$(v_i)$ is a $C$-unconditional basic sequence in $J(h_j)$
for some uniform constant $C$. 
Hence, by passing to the limit above we  get that 
$$
\Vert \sum_{i=1}^m a_iv_i \Vert \leq 3C  \eqno{(4.11)}
$$
for $m \in {\bf N}$.
It follows from (4.11) that $(v_j)$ cannot be a boundedly complete
basis for $[v_i:  i \in {\bf N}]$, 
since $a_{m_{i}} \geq {{\delta}\over 8}$ for $i \in {\bf N}$
by (4.10). Hence [LT,1.c.10] implies that the quasi-reflexive space
$J(h_j)$ must contain an isomorphic copy of $c_0$.
This contradiction completes the proof of Lemma 4.5.

\medskip

\noindent {\sl Claim 2.\  There is $k \in {\bf N}$ and 
$(\alpha _1,\ldots ,\alpha _k) \in (S_{\ell^k_{1}})_{+}$ so that 
the corresponding convex combination  
$A = \sum_{j=1}^k \alpha _jA_j \in J(h_j)^*$ 
given by (4.7) satisfies
$$
\vert A(x)\vert < \delta  \quad {\rm for}\ x \in D. 
$$
}

\noindent {\it Proof of Claim 2.} Let $x \in D$. According to  Step 1 
we may fix a perturbation $x \approx x_0 + \sum_{i=1}^{\infty} a_iy_i$
satisfying conditions (i) - (iii) of Claim 1.
Since  $A_j(x_0) = 0$  by definition, we get that 
$$
\eqalign{\vert A_j(x) \vert \leq \ & \vert  A_j(x_0 + \sum_{i=1}^{\infty} a_iy_i)\vert +
\vert A_j(x - (x_0 + \sum_{i=1}^{\infty} a_iy_i))\vert \cr 
\leq \ & \vert a_jA_j(y_j)\vert 
+ 2\Vert x - (x_0 + \sum_{i=1}^{\infty} a_iy_i)\Vert 
\leq  a_j + {{\delta}\over {2}} \cr
}
$$
for $j \in {\bf N}$. 

Next we use Lemma 4.5 to find $k \in {\bf N}$ and 
$(\alpha _1,\ldots ,\alpha _k) \in (S_{\ell^k_{1}})_{+}$ so that (4.8)
holds. 
Then the above estimate implies that
$$
\vert A(x) \vert \leq \sum_{j=1}^k \alpha _j  \vert A_j(x)\vert 
\leq  \sum_{j=1}^k \alpha _ja_j + 
{{\delta}\over 2}\cdot  \sum_{j=1}^k \alpha _j
< \delta .
$$
This completes the proof of Claim 2.   

\medskip

\noindent {\it Step 3.} 
Fix a positive decreasing sequence $(\varepsilon _i)$ so that
$\sum_{i=1}^{\infty} \varepsilon _i < \varepsilon / 2$. By applying Steps 1 and 2
successively we get a partition ${\bf N} = \bigcup _{j=1}^{\infty} I_j$ 
into successive intervals and
a sequence of functionals $(V_j) \subset J(h_j)^*$
that satisfy the following properties for $j \in {\bf N}$:

\smallskip

\item{(iv)} $\vert V_j(x) \vert < \varepsilon _j/ 2$ for $x \in D$,

\smallskip

\item{(v)} $V_j$ is a convex combination of the type (4.7) 
of averages corresponding to some partition of $I_j$ into successive 
subintervals,

\smallskip

\item{(vi)} $V_j(x) = V_j(\sum_{s\in I_{j}} x(s)e_s)$  
for $x = \sum_{s=1}^{\infty} x(s)e_s \in J(h_i)$.

\smallskip

Write $I_j = [t_j,t_{j+1})$ for $j \in {\bf N}$, where $t_1 = 1$. 
We put  
$$s_0 = e_1,\quad s_1 = \sum_{s=2}^{t_2-1}e_s,\quad 
s_{j} = \sum_{s=t_{j}}^{t_{j+1}-1} e_s\quad {\rm 
for}\ j \geq 2.
$$
Define the linear map $\tilde{V}$ on $J(h_j)$ by
$$
\tilde{V}x = \sum_{j=1}^{\infty} V_j(x)s_{j-1}\quad {\rm for}\ x \in J(h_j). \eqno{(4.12)}
$$
Note that definition (4.12) introduces an additional left "shift" on $J(h_j)$ 
compared to the arguments in sections 2 and 3.  
It is immediate that $V = I - \tilde{V}$  satisfies
$$
\Vert x - Vx\Vert = \Vert \sum_{j=1}^{\infty} V_j(x)s_{j-1}\Vert
\leq 2\sum_{j=1}^{\infty}  \vert V_j(x)\vert < \varepsilon \quad {\rm for\  every}\ x \in D.
$$

We verify below in Claims 3 and 4 that $V = I - \tilde{V}$ 
is the desired weakly compact approximating operator on $J(h_j)$.
The right dominance property (4.2) 
will be essential towards getting a uniform bound for
$\Vert \tilde{V}\Vert$.

\smallskip

\noindent {\sl Claim 3. \quad 
$\Vert V \Vert \leq 5$ (independently of
the subdivisions).
}

\smallskip

\noindent {\it Proof of Claim 3.} 
We estimate $\Vert \tilde{V}\Vert$.
Let $x = \sum_{s=1}^{\infty} x(s)e_s \in J(h_j)$
be finitely supported, and suppose that $supp(x) \subset \bigcup _{j=1}^n I_j$ 
for some $n \in {\bf N}$.
Let  $k_1 < l_1 < \ldots < k_r <  l_r$ be a sequence of coordinates
that norms $\tilde{V}x$ according to (4.1), so that 
$$
\Vert \tilde{V}x\Vert = \Vert \sum_{j=1}^r a_jh_{k_{j}}\Vert ,
$$
where $a_j = \tilde{V}x(k_j) - \tilde{V}x(l_{j})$ for $j = 1,\ldots ,r$.
Since $(h_j)$ is a $1$-unconditional basis of $E$, we may assume 
that $a_j \neq 0$ for $j = 1,\ldots ,r$.

By conditions (v) and (vi) the element
$\tilde{V}x \in J(h_j)$
is  constant on each interval $I_j$ with $j \geq 2$.
In addition, $\tilde{V}x(1) = V_1(x)$ and $\tilde{V}x(s) = V_2(x)$ for 
$2 \leq s < t_2$.
Hence we may assume without loss of generality that no pair
$k_i$ and $l_i$ belongs to the same interval $supp(s_j)$ 
for $i = 1,\ldots ,r$, and that $r \leq  n$.

We claim that
$$
\Vert \sum_{j} a_{2j+1}h_{k_{2j+1}}\Vert \leq 2\Vert x\Vert .
$$
From condition (vi) we get that 
$$
a_{2j+1} =  \tilde{V}x(k_{2j+1}) - \tilde{V}x(l_{2j+1}) = V_{p}(x) - V_{q}(x)
$$
for some $k_{2j+1} \leq p < q$.
By condition (v)  we know that $V_{p}(x)$ is a convex combination of
the coordinates $x(s)$ with $s \in I_{p}$, and similarly for $V_{q}(x)$ 
with respect to $I_{q}$. 
Hence there are coordinates 
$m_{2j+1} \in I_{p}$ and $n_{2j+1} \in I_{q}$, so that
$$
b_{2j+1} \equiv  \vert x(m_{2j+1}) - x(n_{2j+1})\vert  \geq \vert a_{2j+1}\vert 
$$
for each $j$. 
Here $k_1 \leq m_1 < k_3 \leq m_3 < \ldots$, and it is easy to 
convince oneself  that  $m_1 < n_1 < m_3 < n_3 < \ldots$. Thus 
we get from (4.1),  (4.2) 
and the $1$-unconditionality of $(h_j)$ in $E$ that
$$
\Vert \sum_{j} a_{2j+1}h_{k_{2j+1}}\Vert \leq 
\Vert \sum_{j} \vert a_{2j+1}\vert h_{k_{2j+1}}\Vert 
\leq 2 \Vert \sum_{j} b_{2j+1}h_{m_{2j+1}}\Vert \leq  2\Vert x\Vert .
$$
In a similar fashion one has $\Vert \sum_{j} a_{2j}h_{k_{2j}}\Vert \leq 2\Vert x\Vert$,
so that  $\Vert \tilde{V}x\Vert \leq 4\Vert x\Vert$.  
Consequently $\Vert V\Vert \leq 1 + \Vert \tilde{V}\Vert
\leq 5$, which
completes the proof of Claim 3.

\medskip

\noindent {\sl Claim 4. \quad $V = I - \tilde{V} \in W(J(h_j))$.
}

\smallskip

\noindent {\it Proof of Claim 4.}  It suffices to verify that 
$(V(\sum_{j=0}^{n+1} s_j)) = ((I - \tilde{V})(\sum_{j=0}^{n+1} s_j))$ 
is a weak-null sequence in $J(h_j)$, since
$\sum_{j=0}^{n+1} s_j \buildrel{w^*}\over\longrightarrow {\bf 1}$ as 
$n \to \infty$ in $J(h_j)^{**}$. 
Recall for this that $J(h_j)^{**} = J(h_j) \oplus [{\bf 1}]$,
where ${\bf 1} = (1,1,1,\ldots )$,
and that an operator $U \in W(J(h_j))$ if and only if 
$U^{**}({\bf 1}) \in J(h_j)$.

Note first that $\tilde{V}(s_{j+1}) = s_j$ for $j \geq 1$, and that
$\tilde{V}(s_0 + s_1) = s_0$.
It follows that 
$$
(I - \tilde{V})(\sum_{j=0}^{n} s_j) = s_{n+1} \buildrel{w}\over\longrightarrow 0
$$
in $J(h_j)$ as $n \to \infty$. Indeed, if $(s_{n+1})$ is 
not weakly null in $J(h_j)$, then it would contain a subsequence equivalent 
to the unit vector basis of $\ell^1$ by [LT,1.c.9] 
(recall that skipped subsequences
of $(s_{n+1})$ are unconditional, see [BHO,Prop. 2.1.(2)]).
This would contradict the quasi-reflexivity of $J(h_j)$.
The proof of Theorem 4.4 is complete.
$\square$

\medskip

\noindent {\it  Remark 4.6.}  We did not consider the problem whether the dual  
$J(h_j)^*$ always has the W.A.P. in the setting of Theorem 4.4.

\bigskip

\noindent 5. {\bf Permanence properties and further positive results.}

\medskip

In this section we state some simple permanence properties for 
weakly compact approximation,
which imply that certain vector-valued sequence spaces such as 
$\ell^p(\ell^1)$, $\ell^p(J)$ and 
$\ell^1(\ell^p)$ have the W.A.P for $1 < p < \infty$.
These facts will be needed in section 6.

We first recall the following "dual" version of  W.A.P., which has
some applications of its own, see  [T1], [T2].
The Banach space $E$ is said to
have  the {\it inner weakly compact approximation property} (inner W.A.P.
for short)  
if there is a constant $C < \infty$ so that 
$$
\inf \{ \Vert U - UV\Vert : V \in W(E), \Vert V\Vert  \leq C\} = 0 \eqno{(5.1)}
$$
for any  weakly compact operator $U \in W(E,Z)$, where 
$Z$ is an arbitrary Banach space. The inner W.A.P.
is suggested by an analogous operator reformulation of the W.A.P.: 
$E$ has the
W.A.P. if and only if there is $C < \infty$ so that 
$$
\inf \{ \Vert U - VU\Vert : V \in W(E), \Vert V\Vert  \leq C\} = 0 
$$
for any  weakly compact operator $U \in W(Z,E)$, where 
$Z$ is an arbitrary Banach space.  
We recall that the duality between the W.A.P. and the inner
W.A.P. is incomplete.

\medskip

\noindent {\it Examples 5.1.}\ 
{\sl  If $E$ has the inner W.A.P., then $E^*$ has the W.A.P.,
see [T1,3.4]. The converse does not  always hold:
the Johnson-Lindenstrauss space $JL$ fails to have the 
inner W.A.P., but $JL^*$ has the W.A.P., see [T2,1.4].
Moreover, the fact that $E$ has the W.A.P. does not in general imply that  $E^*$ 
has the inner W.A.P. (indeed,
$\ell^1$ has the W.A.P., but 
$\ell^{\infty}$ does not have the inner W.A.P., see [T1,3.5.(ii)]
or Proposition 6.10 below). Note also that $c_0$ has the inner W.A.P.
by [T1,3.5.(ii)], but $c_0$ does not have the W.A.P.}

\medskip

It is a simple fact that the (inner) W.A.P.
is preserved by complementation. We omit the easy arguments.

\medskip

\noindent {\sl Lemma 5.2.\  Suppose that $M \subset E$ is a complemented
subspace, and let $P$ be a projection of $E$ onto $M$.

\smallskip

\noindent (i)\ If $E$ has the W.A.P with constant $C$, then $M$ has the 
W.A.P. with constant $\Vert P\Vert C$.

\smallskip

\noindent (ii)\ If $E$ has the inner W.A.P with constant $C$, then $M$ has the 
inner W.A.P. with constant $\Vert P\Vert C$.
}

\medskip

Let $R$ be a Banach space having a normalized $1$-unconditional Schauder basis 
$(r_j)$ and suppose that $(E_j)$ is a sequence of Banach spaces.
The vector-valued sequence spaces 
(or the $R$-direct sums)
$$
R(E_j) = \{ x = (x_j): x_j \in E_j\ {\rm for}\  j \in {\bf N}, \Vert x\Vert \equiv  
\Vert \sum_{j=1}^{\infty} \Vert x_j\Vert \cdot r_j\Vert _R < \infty\}
$$
provide a suitable setting for our permanence results.
Special cases include the familiar
direct sums $\ell^p(E_j)$ ($1 \leq p < \infty$) and  $c_0(E_j)$.
Let $J_k:  E_k \to R(E_j)$
denote the inclusion map and $P_k$ the natural norm-$1$ projection
of $R(E_j)$ onto  $E_k$ for $k \in {\bf N}$.

\medskip

\noindent {\sl Proposition 5.3.\quad  Let $(E_j)$ be a sequence of Banach spaces
and suppose that $R$ is a reflexive Banach space 
having a normalized $1$-unconditional Schauder basis 
$(r_j)$. Then

\smallskip

\noindent (i) \ $R(E_j)$ has the 
W.A.P. if and only if
$E_j$ has the W.A.P. with a uniform
constant.

\smallskip

\noindent (ii) \  $R(E_j)$ has the 
inner W.A.P. if and only if
$E_j$ has the inner W.A.P. with a uniform
constant. 

\smallskip

\noindent (iii) \  $\ell^1(E_j)$ has the W.A.P. if and only if 
$E_j$ has the W.A.P. with a uniform 
constant.

\smallskip

\noindent (iv) \  $c_0(E_j)$ has the inner W.A.P. if and only 
if $E_j$ has the inner W.A.P. with a uniform 
constant. 
}

\medskip

\noindent {\it Proof.} (i) We put $X = R(E_j)$ for simplicity.
If $X$ has the 
W.A.P.  with constant $C$, then 
the $1$-complemented subspace $E_k \subset X$
has the W.A.P.  with the same constant $C$ for $k \in {\bf N}$
by Lemma 5.2.(i). (The implication "$\Rightarrow$" is checked similarly
for parts (ii), (iii) and (iv).)

Conversely, assume that $E_j$ has the W.A.P.  with 
a uniform constant $C$ for all $j$.
Suppose that $D \subset X$ is a weakly compact subset
and $\varepsilon > 0$. Since $P_kD \subset E_k$ is weakly compact
for $k \in {\bf N}$,  there is by assumption $V_k \in W(E_k)$ so that
$$
\sup_{y \in P_{k}D} \Vert y - V_ky\Vert < {{\varepsilon}\over{2^k}}\quad
{\rm and}\quad  \Vert V_k\Vert \leq C.
$$
Define $V: X \to X$ by $Vx = (V_kx_k)$ for $x = (x_k) \in X$. 
Clearly 
$$
\Vert Vx\Vert = \Vert \sum_{k=1}^{\infty} \Vert V_kx_k\Vert \cdot r_k\Vert _R 
\leq C \Vert \sum_{k=1}^{\infty} \Vert x_k\Vert \cdot r_k\Vert _R = C\Vert x\Vert 
$$
for $x = (x_k) \in X$ by the $1$-unconditionality of $(r_k)$, 
so that $V \in L(X)$ and  $\Vert V\Vert \leq C$.

Let $(x^{(m)}) \subset X$ be a bounded sequence,
where $x^{(m)} = (x_j^{(m)})_{j\in {\bf N}} \in X$ for $m \in {\bf N}$.
Then $x^{(m)}   \buildrel{w}\over\longrightarrow x = (x_j)$ 
in $X$ as $m \to \infty$
if and only if $x_j^{(m)}  \buildrel{w}\over\longrightarrow x_j$ 
as $m \to \infty$ in $E_j$ for $j \in {\bf N}$ (cf.  [L,6.1] for the special
case $\ell^p(E)$). 
The Eberlein-Smulian theorem implies then that
a bounded set $A \subset X$ is relatively weakly compact if and only if
$P_jA$ is relatively weakly compact in $E_j$ for $j \in {\bf N}$.
By applying this fact to $VB_{X}$ we get that $V \in W(X)$,
since $P_j(VB_{X}) = V_jB_{E_{j}}$
is relatively weakly compact for each $j$.
Finally, for $x = (x_k) \in D$ we have
$$
\Vert x - Vx\Vert = \Vert \sum_{k=1}^{\infty} \Vert x_k - V_kx_k\Vert \cdot r_k\Vert _R
\leq \sum_{k=1}^{\infty}{{\varepsilon}\over {2^{k}}} 
= \varepsilon .
$$

\smallskip

\noindent (ii)\  Suppose that $E_j$ has the inner W.A.P. 
with a uniform constant $C$ for all $j$.
Let $U \in W(X,Z)$  and $\varepsilon > 0$ be given, 
where $X = R(E_j)$ and $Z$ is some Banach space.
Put $U_k = UJ_k \in W(E_k,Z)$ for $k \in {\bf N}$. By assumption
there is $V_k \in W(E_k)$ satisfying
$$
\Vert U_k - U_kV_k\Vert < {{\varepsilon}\over{2^k}}\quad
{\rm and}\quad  \Vert V_k\Vert \leq C
$$
for $k \in {\bf N}$. Define $V \in L(X)$
by  $Vx = (V_kx_k)$ for $x = (x_k) \in X$. One verifies as in part (i)
that $\Vert V \Vert \leq C$ and $V \in W(X)$. 

We next estimate $\Vert U - UV \Vert = \Vert U^* - V^*U^*\Vert$.
Observe that $X^* = R^*(E_j^*)$, where the 
biorthogonal sequence $(r_j^*)$ to $(r_j)$ is a $1$-unconditional 
Schauder basis for $R^*$.
Here   $U^* = (U^*_j):  Z^*  \to R^*(E_j^*)$, 
and we get for $z^* \in B_{Z^*}$ that
$$
\Vert U^*z^* - V^*U^*z^*\Vert = 
\Vert \sum_{k=1}^{\infty} \Vert U_k^*z^* - V_k^*U_k^*z^*\Vert \cdot  r_k^*\Vert _{R^*}
\leq \varepsilon  .
$$
Hence $X = R(E_j)$ has the inner W.A.P.

\smallskip

\noindent (iii)  Let $D \subset \ell^1(E_j)$ be a weakly compact set 
and $\varepsilon > 0$.
Using Lemma 5.4 below we fix $n \in {\bf N}$ so that
$$
\sup_{x=(x_{j}) \in D} (\sum_{j=n+1}^{\infty} \Vert x_j\Vert ) 
< {{\varepsilon}\over 2}.
$$
The assumption gives operators $V_j \in W(E_j)$ 
for $j = 1,\ldots ,n$
satisfying $\Vert V_j \Vert \leq C$ and
$\Vert x_j - V_jx_j\Vert < {{\varepsilon}\over {2^{j+1}}}$
for $x = (x_k) \in D$. Define $V \in L(X)$ by
$V(x_k) = (V_1x_1,\ldots ,V_{n}x_{n},0,0,\ldots )$ for
$(x_k) \in X$. Clearly $\Vert V \Vert \leq C$ and $V \in W(X)$.
For $x = (x_j) \in D$ we get that
$$
\Vert x - Vx \Vert \leq \sum_{j=1}^{n} \Vert x_j - V_jx_j\Vert 
+ \sum_{j=n+1}^{\infty} \Vert x_j\Vert \leq \varepsilon .
$$

\smallskip

\noindent (iv)  Put $X = c_0(E_j)$, and  
suppose that $U = (U_k) \in W(X,Z)$
and $\varepsilon > 0$. Here $U_k = UJ_k$ for $k \in {\bf N}$
and $U^* = (U_k^*)$ is a weakly compact operator.
By applying Lemma 5.4 (see below) to the
relatively weakly compact set $U^*(B_{Z^{*}}) \subset X^* = \ell^1(E_j^*)$
there is $n \in {\bf N}$ so that
$$
\sup_{z^* \in B_{Z^{*}}} (\sum_{j=n+1}^{\infty} \Vert U_j^*z^* \Vert ) 
< {{\varepsilon}\over 2}.
$$
Fix $V_j \in W(E_j)$ so that $\Vert U_j - U_jV_j\Vert < {{\varepsilon}\over {2^{j+1}}}$ 
and $\Vert V_j\Vert \leq C$ for $j = 1,\ldots ,n$.
Define again $V \in W(X)$ by 
$V(x_k) = (V_1x_1,\ldots ,V_{n}x_{n},0,0,\ldots )$ for
$(x_k) \in X$, so that $\Vert V \Vert \leq C$.
For $z^* \in B_{Z^{*}}$ we get that
$$
\Vert U^*z^* - V^*U^*z^*\Vert \leq 
\sum_{j=1}^{n} \Vert U_j^*z^* - V_j^*U_j^*z^* \Vert 
+ \sum_{j=n+1}^{\infty} \Vert U_j^*z^*\Vert \leq  
\sum_{j=1}^{n} \Vert U_j - U_jV_j \Vert + {{\varepsilon}\over 2} \leq \varepsilon .
$$
Thus $\Vert U - UV\Vert \leq \varepsilon$, so that $c_0(E_j)$ has the 
inner W.A.P.
$\square$

\medskip

The following auxiliary fact was used 
in the proof of parts (iii) and (iv) of Proposition 5.3.
We sketch the argument of this well known result
for completeness.

\medskip

\noindent {\sl Lemma 5.4.} {\sl Let $(E_j)$ be a sequence of Banach spaces
and suppose that $D \subset \ell^1(E_j)$ is 
a weakly compact set. 
Then for any
$\delta > 0$ there is $n = n(D,\delta ) \in {\bf N}$ so that
$$
\sup_{y=(y_{j}) \in D} (\sum_{j=n}^{\infty} \Vert y_j\Vert ) 
\leq \varepsilon . \eqno{(5.2)}
$$
}

\noindent {\it Proof.}  Suppose to the contrary that (5.2) does not hold:
there is a $\delta > 0$ so that for each $n \in {\bf N}$
there is some  $y = (y_j) \in D$ satisfying 
$\sum_{j=n}^{\infty} \Vert y_j\Vert > \delta$.
Hence there are sequences
$(p_n) \subset {\bf N}$ and  $(y^{(n)}) \subset D$
so that
$\sum_{j=p_{n}+1}^{p_{n+1}} \Vert y_j^{(n)}\Vert > \delta$ 
for $n \in {\bf N}$. Here $y^{(n)} = (y^{(n)}_j)_{j\in {\bf N}} \in \ell^1(E_j)$ for 
$n \in {\bf N}$.
It is then easy to verify that a subsequence 
of $(y^{(n)})$ is equivalent to the unit vector basis of $\ell^1$
(this will contradict the weak compactness of $D$). 
$\square$

\medskip

The following novel examples of
concrete spaces that have the W.A.P. (or
the inner W.A.P.)  are immediate from Proposition 5.3.

\medskip

\noindent {\sl Corollary 5.5.} {\sl  Let $1 < p < \infty$. Then

\smallskip

\noindent (i)\quad  $\ell^p(\ell^1)$ and $\ell^1(\ell^p)$ have the W.A.P.,

\noindent (ii)\quad $\ell^p(c_0)$ and $c_0(\ell^p)$ have the inner W.A.P.,

\noindent (iii)\quad  $\ell^p(J)$ and $\ell^p(J^*)$ have the W.A.P.
}

\medskip

\noindent {\sl Proof.}  Recall that $\ell^1$ has the W.A.P., and that 
$c_0$ has the inner W.A.P.
(see [T1,3.5.(ii)]). Part (iii) follows from Theorems 2.2, 3.3
and Proposition 5.3.(i).
$\square$

\medskip

\noindent {\sl Remark 5.6.}\  The fact that $\ell^2(J)$ has the W.A.P. 
sheds some further light on a result of [GST].
Let $E$ be a Banach space and define the "residual" operator
$R(S) \in L(E^{**}/E)$  by 
$$
R(S)(x^{**} + E) = S^{**}x^{**} + E\quad {\rm for}\ x^{**} \in E^{**},\ 
S \in L(E).
$$
It is known (cf. [GST,1.4]) that $\Vert R(S)\Vert \leq \omega (S) \leq \Vert S \Vert _w$
for $S \in L(E)$, where $\Vert S \Vert _w \equiv {\rm dist}(S,W(E))$ 
and $\omega (\cdot )$ is the measure of weak non-compactness (cf. the Introduction).

According to [GST,2.6]  there is a sequence $(S_n) \subset L(\ell^2(J))$
so that $\Vert S_n  \Vert _w = 1$ for all $n$, but $\Vert R(S_n)\Vert
\to 0$ as $n \to \infty$. The precise relation between $\omega$ and 
$\Vert \cdot \Vert _w$ on $L(\ell^2(J))$ was not resolved in [GST]. 
Now an inspection of the arguments of 
Proposition 5.3.(i) and [AT,Thm. 1]  reveals that in fact 
$$
\omega (S) \leq \Vert S \Vert _w  \leq 2\  \omega (S), \quad S \in L(\ell^2(J)).
$$

\medskip

Another natural permanence problem, which we did not pursue here,
concerns the W.A.P for the Bochner spaces $L^p(E)$.  

\medskip

\noindent {\sl Problem 5.7.\quad  Does 
$L^p(E)  = L^p([0,1];E)$ have the  W.A.P (resp., the inner W.A.P.)
whenever $E$ has the W.A.P. (resp., the inner W.A.P.)
and  $1 < p < \infty$ ? (The cases $p = 1$ and $p = \infty$ are 
excluded by known facts, cf. Proposition 6.10).
}

\bigskip

\noindent 6.\quad {\bf James'  tree space $JT$ does not have the W.A.P.
and related examples.}

\medskip

This section provides concrete answers to
various natural questions about the class of spaces having the 
W.A.P.  We first recall a couple of notions. 
Let $1 \leq p < \infty$ be fixed. The Banach space $E$ is {\it $\ell^p$-saturated} if
every  infinite-dimensional subspace $M \subset E$ 
contains an isomorphic copy of $\ell^p$. The space $E$ is 
{\it somewhat reflexive}
if every infinite-dimensional subspace $M \subset E$ contains a
reflexive infinite-dimensional subspace. 
(Here "subspace" always means a closed linear subspace.)

The quasi-reflexive H.I. space $E$ from  [ArT,Prop. 14.10] (cf. Example 4.1)
that fails the W.A.P.
yields a striking counterexample to the following question stated in [AT].

\medskip

\noindent {\it Question 6.1.\  } [AT,p. 370] 
{\it Suppose that the quotient $E^{**}/E$ is
reflexive. Does $E$ have the W.A.P.?}

\medskip

Since $\ell^2(J)$ has the W.A.P. by Corollary 5.5.(iii),
and $\ell^2(J)^{**}/\ell^2(J) = \ell^2(J^{**}/J) =  \ell^2$, 
there are spaces $E$ with  the W.A.P. and $E^{**}/E$ 
reflexive and infinite-dimensional.
Theorem 6.5 below yields a concrete space 
$Y$ without the W.A.P.  for which $Y^{**}/Y = \ell^2$. 

Our next question addresses another potential
extension of the fact that reflexive spaces have the W.AP.

\medskip

\noindent {\it Question 6.2.\ 
Suppose that $E$ is a somewhat reflexive space 
that has the bounded approximation property (B.A.P.).
Does $E$ have the W.A.P.?
}  

\medskip

The answer to Question 6.2 can be deduced from known results 
(Theorem 6.5 below contains a different, 
$\ell^2$-saturated example).

\medskip

\noindent {\sl Example 6.3.\ 
Let $E$ be the separable, somewhat reflexive ${\cal L}^{\infty}$-space
constructed by Bourgain and Delbaen, see [B,Ch. III] for a description.
Then $E$ has the B.A.P., 
but $E$ does not have the W.A.P. (see 
Proposition 6.10 below or [AT,Cor. 3]). 
} 
$\square$

\medskip

James'  tree space $JT$ was introduced by James [J2] 
as a useful variation of the ideas underlying $J$, 
and its properties were further analyzed e.g. by
Lindenstrauss and Stegall [LS].  There is a systematic exposition 
of the properties of $JT$ in [FG, chapter 3].
The fact  that $J$ has the W.A.P. (section 2) suggests the following problem.

\medskip

\noindent {\it Question 6.4.\  Does $JT$ have the W.A.P.?
}

\medskip

The main result of this section (Theorem 6.5) establishes that 
$JT$ does not have the W.A.P., where  
$JT^{**}/JT$ is a non-separable Hilbert space.
We recall the definition of $JT$ and fix some relevant notation. 
Let 
$$
{\cal T} = \bigcup _{n=0}^{\infty} \{0,1\}^n
$$ 
be the infinite binary tree  equipped with the natural partial order. The {\it nodes}
$\alpha \in {\cal T}$ satisfy
$\alpha = \emptyset$ or $\alpha = (\alpha _1,\ldots ,\alpha _n)$ for some
$n \in {\bf N}$, where $\alpha _j = 0$ or $\alpha _j = 1$ for $j \in \{1,\ldots ,n\}$.
The {\it length} $\vert \alpha \vert$ of $\alpha = 
(\alpha _1,\ldots ,\alpha _n) \in {\cal T}$ is $n$.
Given $\alpha \in {\cal T}$ let $f_{\alpha}: {\cal T} \to {\bf R}$ 
be defined by $f_{\alpha}(\alpha ) = 1$ 
and $f_{\alpha}(\beta ) = 0$ for  $\beta \neq \alpha$. James'  tree space $JT$
consists of the functions $\sum_{\alpha} a_{\alpha} f_{\alpha}:  {\cal T} \to {\bf R}$
for which the norm
$$
\Vert \sum_{\alpha} a_{\alpha} f_{\alpha}\Vert = \sup_{k; S_{1},\ldots ,S_{k}} 
(\sum_{j=1}^k S_j^*(\sum_{\alpha} a_{\alpha} f_{\alpha})^2)^{1/2} < \infty , \eqno{(6.1)}
$$
where the supremum is taken over disjoint segments $S_1,\ldots ,S_k$
of ${\cal T}$ and $k \in {\bf N}$. A {\it segment} $S \subset {\cal T}$ has the form
$S = \{\gamma \in {\cal T}: \alpha \leq \gamma \leq \beta \}$ 
for given $\alpha , \beta \in  {\cal T}$ with $\alpha \leq \beta$, and 
$S^*(\sum_{\alpha} a_{\alpha} f_{\alpha}) =
\sum_{\alpha \in S} a_{\alpha}$ for $\sum_{\alpha} a_{\alpha} f_{\alpha} \in JT$. 
It is known that  $(f_{\alpha})_{\alpha \in {\cal T}}$ is a monotone boundedly complete
basis for $JT$ (ordered by increasing length of the nodes 
and from "left to right").

A {\it branch} $B \subset {\cal T}$ is a maximal infinite order interval
starting at $\emptyset$.  
A branch $B$ determines the norm-$1$ functional 
$S_B^* \in JT^*$ defined by $S_B^*(\sum_{\alpha} a_{\alpha} f_{\alpha})
= \sum_{\alpha \in B} a_{\alpha} f_{\alpha}$. Let 
$\Gamma$ be the uncountable collection of all branches of
${\cal T}$.  Then
$JT^{**}/JT = \ell^2(\Gamma )$ isometrically, and
$$
JT^* = [\{f_{\alpha}^*: \alpha \in {\cal T}\} \cup \{S_B^*: 
B \in \Gamma ],
\eqno{(6.2)}
$$
see  [LS,Thm. 1] or [FG,3.c.3].
Here $(f_{\alpha}^*)_{\alpha \in {\cal T}}$ is the biorthogonal sequence to
$(f_{\alpha})_{\alpha \in {\cal T}}$  in $JT^*$.
Recall further that $JT$ is $\ell^2$-saturated, see [J2,Thm.] or [FG,3.a.8].

The following  special notation will be  convenient. 
If $\alpha = (\alpha_1,\ldots ,\alpha _n)\in {\cal T}$,
then $\alpha 0 = (\alpha_1,\ldots ,\alpha _n,0)$ is the {\it left successor} 
(or left concatenation) and 
$\alpha 1 = (\alpha_1,\ldots ,\alpha _n,1)$ the {\it right successor} of $\alpha$.
For any $n \in {\bf N}$ we put
$$
{\cal T}_n = \{ \alpha \in {\cal T}: {\rm there\ are\ at\ most}\ n\  1's\  {\rm in}\  \alpha\}.
$$
Thus  $\alpha \in {\cal T}_n$ if the node $\alpha$ contains at most $n$ "right turns".
Put  $X_n = [f_{\alpha}: \alpha \in {\cal T}_n]$ for $n \in {\bf N}$. 
Note that $X_n \subset JT$ is a 
$1$-complemented subspace, where the restriction $x \mapsto x_{|{\cal T}_{n}}$
defines the natural projection onto $X_n$. Indeed, if $S_1,\ldots ,S_k \subset {\cal T}$
are disjoint segments, then $S_j \cap {\cal T}_n$ are disjoint segments 
(possibly empty) in ${\cal T}_n$ for $j = 1,\ldots ,k$.

\smallskip

We are ready for our main results about $JT$. 
Below parts (i) and (ii) together imply  that $X_n$ 
has the W.A.P. for all $n \in {\bf N}$,
but where the smallest constant $C$ in (1.1) is proportional to $\sqrt{n}$. 
Parts (iii) and (iv) are only based on (ii), but (i) will  become useful later
(see Example 6.8 and Remark 6.9).
The space $Y$ in part (iv) is an example,
where $Y$ fails to have the W.A.P. and the reflexive quotient $Y^{**}/Y$ 
is much "smaller" than $JT^{**}/JT$.  The quasi-reflexive space 
from [ArT,Prop. 14.10]  yields an optimal negative answer to Question 6.1
in terms of minimizing $dim(Y^{**}/Y)$, 
but our examples are easier.
They are further witnesses 
that the quotient  $E^{**}/E$ alone 
does not decide the W.A.P. of $E$.

\medskip

\noindent {\it Theorem 6.5. Let $X_n = [f_{\alpha}: \alpha \in {\cal T}_n] \subset JT$ 
be as above. Then the following
properties hold.

\smallskip

\noindent (i)\ $X_n$ has the W.A.P. with constant (at most) $3\sqrt{n}$ for $n \in {\bf N}$.

\smallskip

\noindent (ii)\  There is a uniform constant $C > 0$ with the following property:
for each $n \in {\bf N}$ there is a weakly compact set $D_n \subset X_n$
so that if 
$$
\sup_{x \in D_{n}} \Vert x - Vx\Vert  < {1\over {10}} \quad  {\rm and}\ V \in W(X_n), 
\eqno{(6.3)}
$$
then $\Vert V\Vert \geq C\sqrt{n}$.

\smallskip

\noindent (iii)\ $JT$ does not have the W.A.P.

\smallskip

\noindent (iv)\  If $Y = \ell^2(X_n)$, then $Y^{**}/Y = \ell^2$ isometrically and
$Y$ does not have the W.A.P.
}

\medskip

\noindent {\it Proof.}  We say that  $B\subset {\cal T}$
is a branch {\it starting at the node} $\alpha \in {\cal T}$, if $B$ is a 
maximal infinite linearly ordered set so that $\gamma \geq \alpha$ for
all $\gamma \in B$.
It is convenient to fix, for each $n \in {\bf N}$, a partition
$$
{\cal T}_n = \bigcup_{j=1}^{\infty} B_j^{(n)}
$$ 
into disjoint branches, where every
$B_j^{(n)} = (\alpha^{(j,n)} ,\alpha^{(j,n)}0,\alpha^{(j,n)}00,\ldots )$ 
is the "always left" branch
starting at some node 
$\alpha^{(j,n)} = (\alpha^{(j,n)} _1,\ldots ,\alpha^{(j,n)} _k) \in {\cal T}_n$ 
with $\alpha^{(j,n)} _k = 1$. We may enumerate these branches by requiring that
$B_{r+1}^{(n)}$ starts at the first node $\alpha \in {\cal T}_n \setminus
\bigcup _{j=1}^r B_j^{(n)}$
(enumerated according to increasing length and 
from "left to right").

\medskip

\noindent (i)   We will apply Proposition 5.3.
Fix $n \in {\bf N}$ and put
$B_j \equiv  B_j^{(n)}$ for notational simplicity as $j \in {\bf N}$. 
Note that 
$Y_j \equiv  [f_{\alpha}: \alpha \in  B_j] = J$ isometrically for $j \in {\bf N}$
according to (6.1) and (3.1).
Let $P_j \in L(X_n)$ be the natural norm-$1$ projection
onto $Y_j$ corresponding to the restriction $x \mapsto x_{| B_{j}}$ for 
$j \in {\bf N}$.

Define a linear map $T: X_n \to \ell^2(J)$ by $Tx = (P_jx)$ for $x \in X_n$.
It will be enough to show that $T$ is an isomorphism 
satisfying $\Vert T \Vert \cdot \Vert T^{-1}\Vert \leq \sqrt{n}$.
Indeed, recall that according to Theorem 2.2 and Proposition 5.3.(i) 
the direct sum $\ell^2(J)$ has the W.A.P. with constant $C \leq 3$  
as defined by (1.1).
It is then straightforward to check (using the isomorphism $T$) that 
$X_n$ has the W.A.P.  with some constant $\tilde{C} \leq 3\sqrt{n}$.

We claim that the following estimates hold, where 
the left-hand inequality of (6.4)
states that $T$ is well-defined $X_n \to \ell^2(J)$.

\smallskip

\noindent {\it Claim 1.  If $x \in X_n$, then
$$
(\sum_{j=1}^{\infty} \Vert P_jx\Vert ^2)^{1/2} \leq 
\Vert x \Vert \leq \sqrt{n} \cdot (\sum_{j=1}^{\infty} \Vert P_jx\Vert ^2)^{1/2}. 
\eqno{(6.4)}
$$
}

\smallskip

\noindent {\it Proof of Claim 1.}  We may write $x \in X_n$ coordinatewise as
$x = \sum_{j=1}^{\infty} P_jx =  \sum_{j=1}^{\infty} x_{| B_{j}}$,
since the branches $\{B_j:  j \in {\bf N}\}$ form a partition of 
${\cal T}_n$. The left-hand inequality in (6.4) is then obvious by selecting
segments that approximately norm each $P_jx$ and which are 
wholly contained in $B_j$.

Suppose next that $S_1,\ldots ,S_m \subset {\cal T}_n$ are disjoint segments.
According to (6.1) we must show that
$$
\sum_{j=1}^m S_j^*(x)^2 \leq n\cdot \sum_{r=1}^{\infty} \Vert P_rx\Vert ^2
$$
for $x \in X_n$.
Recall that the nodes $\alpha \in {\cal T}_n$ have at most $n$ right turns,
so that 
$$
n(j) \equiv \vert \{r \in {\bf N}: S_j \cap B_r \neq \emptyset\}\vert  \leq n
$$
for $j = 1,\dots ,m$. 
Write the resulting intersected segments  as 
$S_{j,j(1)},\ldots ,S_{j,j(n)}$, where we put $S_{j,j(r)} = \emptyset$ if $n(j) < r \leq n$. 
We may thus write
$S_j^*(x) = \sum_{r=1}^n S_{j,j(r)}^*(x)$ for $x \in X_n$ and $j = 1,\ldots ,m$
(observing the convention that $S_{j,j(r)}^*(x) = 0$ if  $S_{j,j(r)} = \emptyset$).
Here  $(\sum_{r=1}^n S_{j,j(r)}^*(x))^2 \leq n\cdot (\sum_{r=1}^n S_{j,j(r)}^*(x)^2)$
by H\"older's inequality. Hence we get from the above that
$$
\sum_{j=1}^m S_j^*(x)^2 = \sum_{j=1}^m ( \sum_{r=1}^n S_{j,j(r)}^*(x))^2
\leq n\cdot \sum_{j=1}^m(\sum_{r=1}^n S_{j,j(r)}^*(x)^2)
\leq n\cdot  \sum_{s=1}^{\infty} \Vert P_sx\Vert ^2.
$$
For the  last estimate regroup the finite sum into those 
of  the disjoint segments 
$\{S_{j,j(r)}: j = 1,\ldots ,m,\ r = 1,\ldots ,n\}$
that  lie inside any given branch $B_s$ for $s \in {\bf N}$.

\medskip

\noindent (ii)  Let $n \geq 6$ be fixed. For simplicity we put again 
$B_j \equiv B_j^{(n)}, j \in {\bf N}$, 
for the partition ${\cal T}_n = \bigcup_{j=1}^{\infty} B_j^{(n)}$
that was fixed at the beginning of the proof.
Consider the subset
$$
D_n = \bigcup _{\alpha \in {\cal T}_{n}} \{f_{\alpha} - f_{\alpha 0}, 
f_{\alpha 0} - f_{\alpha 00},\ldots \}
= \bigcup _{j=1}^{\infty} \{ f_{\alpha} - f_{\alpha 0}: \alpha \in B_j\} \subset X_n.
$$
Here the sequences $( f_{\alpha} - f_{\alpha 0})_{\alpha \in B_{j}}$ 
are formed by the consecutive differences
along the "always left" branches $B_j$ in ${\cal T}_n$ for $j \in {\bf N}$.

\medskip

\noindent {\it Claim 2.\quad The set $D_n \cup \{0\}$ is weakly compact in $X_n$.
}

\medskip

\noindent {\it Proof of Claim 2.} We will verify that 
any sequence $(x_m) = (f_{\alpha _{m}} - f_{\alpha _{m}0})$
of distinct points from $D_n$ contains a weak-null subsequence
$(x_{m_{k}})$ in $JT$.  

There is no loss of generality to assume,
by applying Ramsey's  classical theorem and
passing to a subsequence of $(x_m)$, that either

\smallskip

\item{(6.5)} all the nodes $\alpha _m$ lie on a single branch  
$B$ of ${\cal T}$, or

\smallskip

\item{(6.6)}  all the nodes $\alpha _m$ are pairwise incomparable.

\smallskip

If (6.5) holds, then $(x_m)$ is equivalent to a subsequence of 
the shrinking basis $(e_n)$ of $J$ given by (2.1), and hence it is weakly
null (cf. [FG,2.c.10]). 
If (6.6) holds, then $(x_m)$ is equivalent to the unit vector basis
of $\ell^2$, and hence it is again weakly null. 
Thus Claim 2 holds.

\smallskip

Suppose next that  $V \in L(X_n)$ is a weakly
compact operator satisfying (6.3) for the weakly compact set $D_n \cup \{0\}$.
For simplicity, let  $(f_m)$ stand for the node basis of a given
branch $B_r$ of the partition of ${\cal T}_n$.
We make a preliminary observation.

\smallskip

\noindent {\it Fact.} Given $\delta > 0$ there is a 
sequence of disjointly supported
convex blocks $(g_j)$ of $(f_m)$ so that 
$$
\Vert Vg_j - Vg_i \Vert < \delta ,\quad i \neq j. \eqno{(6.7)}
$$
Here  $\Vert g_j \Vert = 1$ for $j \in {\bf N}$ in view of (6.1).

\smallskip

\noindent Indeed, the weak compactness of $V$ yields
a subsequence $(f_{m_{i}})$ so that 
$Vf_{m_{i}} \buildrel{w}\over\longrightarrow x \in JT$ as $i \to \infty$.
Then  Mazur's theorem gives
a sequence of disjointly supported
convex blocks $(g_j)$ of $(f_{m_{i}})$ so that 
$\Vert Vg_j - Vg_i \Vert < \delta$ whenever  $i \neq j$.

\smallskip

Let $[t]$ denote the integer part of $t > 0$.
We successively apply the preceding Fact to 
$[n/2]$ "adjacent" branches in ${\cal T}_n$, in the manner described below, 
to get the  element 
$$
x_n = \sum_{j=1}^{[n/2]} (g_{2j} - g_{2j-1}) +  
\sum_{j=1}^{[n/2]} (f_{\alpha_{j}} - f_{\alpha_{j}0}) \in X_n.   \eqno{(6.8)}
$$
The differences $g_{2j} - g_{2j-1}$ and $f_{\alpha_{j}} - f_{\alpha_{j}0}$
are successively  chosen as follows for $j = 1,\ldots ,[n/2]$:

\smallskip

\item{(6.9)} $g_1$ and $g_2$ are normalized convex blocks
on the node basis determined by the "always left" branch $B_1$, 
their supports satisfy 
$\max supp(g_1) < \gamma _1 < \min supp(g_2)$ 
for some node $\gamma _1 \in B_1$,
and $\Vert Vg_{2} - Vg_{1}\Vert < {1\over {10}}$. 
(Here the support of the convex combinations
is with respect to the node basis).

\smallskip

\item{(6.10)} $\alpha _1 = \gamma_11$ (the right successor of $\gamma _1$).

\smallskip

To continue repeat the above procedure by applying 
(6.7) to the "always left" branch in ${\cal T}_n$
starting from the node $\alpha _11 = \gamma_111$ (the right successor of $\alpha _1$).
A picture will be helpful at this stage.
This construction can be performed $[n/2]$ times, since 
${\cal T}_n$ allows at most $n$ right turns.

\smallskip

\noindent {\it Claim 3.} 
$$
\Vert x_n\Vert \leq \sqrt{6n} \quad {\rm and}\quad  \Vert Vx_n \Vert \geq {n\over 3}.
\eqno{(6.11)}
$$
Clearly (6.11) yields that
$\Vert V \Vert \geq {1\over {\sqrt{6n}}} \cdot {n \over 3} = 
{1\over {3\sqrt{6}}} \cdot \sqrt{n}$,
which gives (ii) with $C = {1\over {3\sqrt{6}}}$.

\smallskip

\noindent {\it Proof of  Claim 3.} Let $S_1,\ldots ,S_m \subset {\cal T}_n$ be given disjoint 
segments.  We have to verify that
$$
(\sum_{j=1}^m S_j^*(x_n)^2)^{1/2} \leq \sqrt{6n}.
$$
Note that $x_n$ is a sum of $4[{n\over 2}]$ normalized blocks in $JT$,
namely the convex blocks $g_i$ for $i = 1,\ldots ,2[{n\over 2}]$, and 
$f_{\alpha _{j}}$ and $f_{\alpha _{j}0}$ for $j = 1,\ldots ,[{n\over 2}]$.
From the iterative construction
of $x_n$ it follows that for each segment $S_j$ there are at most
$3$ non-empty disjoint segments $S_{j,1}, S_{j,2}, S_{j,3} \subset S_j$
so that

\smallskip

\item{(6.12)} $S_{j}^*(x_n) = \sum_{i=1}^3 S_{j,i}^*(x_n)$,

\smallskip

\item{(6.13)} each $S_{j,i}$ is contained in the smallest segment 
containing one of the blocks forming $x_n$.

\smallskip

Let $T_s$ be the smallest segment in ${\cal T}_n$ containing $supp(g_s)$
for $s = 1,\ldots ,2[{n\over 2}]$. Note that 
$\sum_{r=1}^p U_r^*(x_n)^2 
\leq \Vert g_s\Vert ^2 = 1$ for each $s$, whenever $U_1,\ldots ,U_p$ are
disjoint segments contained in $T_s$.
Hence it follows from H\"older's inequality and (6.12),(6.13) that
$$
(\sum_{j=1}^m S_j^*(x_n)^2)^{1/2} \leq \sqrt{3} \sqrt{4[{n\over 2}]}
\leq \sqrt{6n}.
$$
This yields the first estimate in (6.11).

Note for the second estimate in (6.11) that according 
to assumption (6.3) one has
$$
V(f_{\alpha_{j}} - f_{\alpha_{j}0})
= f_{\alpha_{j}} - f_{\alpha_{j}0} + z_j,
$$
where $\Vert z_j\Vert  <  {1\over {10}}$
for $j = 1,\ldots ,[n/2]$, since  
$f_{\alpha_{j}} - f_{\alpha_{j}0} \in D_n$.
Let $S \subset {\cal T}$ be a segment so that $\alpha _j \in S$, but
its left successor $\alpha _j0 \notin S$ for $j = 1,\ldots ,[n/2]$. Then we
get that
$$
\Vert \sum_{j=1}^{[n/2]}(f_{\alpha_{j}} - f_{\alpha_{j}0})\Vert \geq 
\vert S^*(\sum_{j=1}^{[n/2]}(f_{\alpha_{j}} - f_{\alpha_{j}0}))\vert 
= [{n\over 2}].
$$
Since $\Vert V(g_{2j} - g_{2j-1})\Vert  < {1\over {10}}$ for $j = 1,\ldots ,[n/2]$
by construction, we obtain that
$$
\eqalign{\Vert Vx_n \Vert = \ & \Vert \sum_{j=1}^{[n/2]} V(g_{2j} - g_{2j-1}) +  
\sum_{j=1}^{[n/2]} V(f_{\alpha_{j}} - f_{\alpha_{j}0})\Vert \cr 
\geq \ &  \Vert \sum_{j=1}^{[n/2]}(f_{\alpha_{j}} - f_{\alpha_{j}0})\Vert -
\sum_{j=1}^{[n/2]} \Vert z_j\Vert - \sum_{j=1}^{[n/2]}\Vert  V(g_{2j} - g_{2j-1})\Vert \cr
\geq \ & [{n\over 2}] - {n\over 2}({1\over {10}} + {1\over {10}}) 
\geq {n\over 3}.
}
$$

\medskip

\noindent (iii)\  This fact follows from part (ii) and Lemma 5.2.(i), 
since $D_n \subset X_n \subset JT$,
where $X_n$ is $1$-complemented in $JT$ 
for all $n \in {\bf N}$

\medskip

\noindent (iv)\ The direct sum $Y = \ell^2(X_n)$ does not have the W.A.P.  in view
of Proposition 5.3.(i), since according to part (ii) the spaces 
$X_n$ do not have the W.A.P. with 
a uniform constant. 
A modification of the corresponding argument for
$JT$ in [LS,Thm. 1] (see also [FG,3.c.3]) yields that $X_n^{**}/X_n = \ell^2$
isometrically for all $n \in {\bf N}$. Hence 
$\ell^2(X_n)^{**}/\ell^2(X_n) = \ell^2(X_n^{**}/X_n) = \ell^2$.
$\square$

\medskip

\noindent {\sl Remark 6.6.}\  Lindenstrauss and Stegall [LS]
defined a function space version of
$J$. James'  function space $JF$ does not have the W.A.P., 
since the separable space $JF$ 
contains a (complemented) copy of $c_0$, see [LS,p. 95].

\medskip

The property defined by (1.1) should more precisely be called
the {\it bounded} W.A.P.  We say that $E$ has the {\it unbounded W.A.P.}
if for every weakly compact set $D \subset E$ and $\varepsilon > 0$ 
there is $V \in W(E)$ satisfying
$$
\sup_{x \in D} \Vert x - Vx\Vert < \varepsilon . \eqno{(6.14)}
$$
The reason for our unorthodox terminology is that
the known applications of weakly compact approximation 
are related to the property defined by (1.1), rather than the one by (6.14).
Recall that there are spaces that have the (finite rank) approximation
property A.P., but not the B.A.P., see [LT,1.e].
This raises another problem.

\medskip

\noindent {\sl Question 6.7.\quad  Is there a space $E$ that has 
the unbounded W.A.P., but not the W.A.P.?
}

\medskip

It turns out that Theorem 6.5 yields concrete examples of this kind, 
so that the unbounded W.A.P. is a strictly weaker notion than the W.A.P.

\medskip

\noindent {\it Example 6.8.\  Let $X_n = [f_{\alpha}: \alpha \in {\cal T}_n]
\subset JT$ be the spaces from  Theorem 6.5 for $n \in {\bf N}$, and
let $Z = \ell^1(X_n)$ be their direct $\ell^1$-sum. Then 
$Z$ has the unbounded W.A.P., but not the W.A.P.
}

\medskip

\noindent {\it Proof.} Proposition 5.3.(iii) yields that $Z = \ell^1(X_n)$
does not have W.A.P., since the spaces $X_n$ do not have the 
W.A.P. with a uniform constant according to Theorem 6.5.(ii).
On the other hand, since $X_n$ has the W.A.P. for all $n \in {\bf N}$
by Theorem 6.5.(i), a simple modification of the argument for Proposition
5.3.(iii) implies that $Z = \ell^1(X_n)$ does have the unbounded W.A.P.
Indeed, recall that the relevant approximating operators $V \in W(Z)$ were 
defined by $Vx = (V_1x_1,\ldots ,V_nx_n,0,0,\ldots )$, 
$x = (x_k) \in Z$,  for suitably
chosen $n \in {\bf N}$ and $V_j \in W(X_j)$ for $j = 1,\ldots ,n$.
$\square$

\medskip

\noindent {\it Remark 6.9.}\  The space 
$JT$ does not even have the unbounded W.A.P.  Indeed, let
$$
\tilde{D} = \bigcup _{n=1}^{\infty} D_n \cup \{0\} \subset  \ell^2(X_n)
$$
be the coordinatewise union in the direct $\ell^2$-sum, 
where  the weakly compact sets  
$D_n \subset X_n$ are
those of the proof of Theorem 6.5.(ii) for $n \in {\bf N}$. 
The set $\tilde{D}$ is relatively weakly compact in $\ell^2(X_n)$
(cf. the proof of Proposition 5.3.(i)).
Note that  $\ell^2(X_n) \subset \ell^2(JT)$, where $\ell^2(JT)$ embeds
as a complemented subspace of $JT$, see [FG,3.a.17]. Fix a 
linear embedding $T: \ell^2(JT) \to JT$, and 
a projection $P$ of $JT$ onto $T(\ell^2(X_n))$.

Suppose  that for any $\varepsilon > 0$
there is $V \in W(JT)$ satisfying $\Vert x - Vx\Vert < \varepsilon$
for all $x \in T(\tilde{D})$. 
It is then easy to check that for every $n \in {\bf N}$ 
there is $V_n \in W(X_n)$, so that 
$$
\Vert V_n\Vert \leq C\Vert V \Vert \quad {\rm and} \quad  
\sup_{z \in D_{n}} \Vert z - V_nz\Vert < c \cdot \varepsilon ,
$$ 
where $C > 0$ and $c > 0$ are uniform constants 
that only depend
on $\Vert T\Vert$, $\Vert T^{-1}\Vert$ and $\Vert P \Vert$.
This contradicts Theorem 6.5.(ii) with  
$\epsilon > 0$ small enough and $n \in {\bf N}$
large enough. 
$\square$

\medskip

We next state for completeness 
two simple conditions which guarantee that spaces with the 
Dunford-Pettis property fail  to have the (inner) W.A.P.
(see also [AT, Prop. 2] and [T1,3.3]). Recall that a Banach space $E$ has the 
{\it Dunford-Pettis property} (DPP) if $\Vert Vx_n\Vert \to 0$ as $n \to \infty$
whenever $V \in W(E,F)$ and $(x_n) \subset E$ is a weak-null sequence.
The space $E$ has the {\it Schur property} if 
$\Vert x_n\Vert \to 0$ as $n \to \infty$ for every weak-null 
sequence $(x_n) \subset E$.
The survey [Di] contains a lot of information about 
the Dunford-Pettis and the Schur properties.
The known facts that $L^1(0,1)$ and $L^{\infty}(0,1)$ have neither the
W.A.P. nor the inner W.A.P., as well as many additional examples, 
can be recovered from the following proposition.

\medskip

\noindent {\it Proposition 6.10.} {\sl
Let $E$ be a Banach space having the DPP.

\smallskip

\noindent (i) If $E$ has the W.A.P., then $E$ has the Schur property.
In particular, if  $E$ contains an infinite-dimensional reflexive subspace
$M$, then $E$ fails to have the W.A.P.

\smallskip

\noindent (ii) If $E$ has an infinite-dimensional reflexive quotient space
$E/M$, then $E$ fails to have the inner W.A.P.
}

\medskip

\noindent {\it Proof.} (i)\ If $E$ does not have the Schur property, then there
is a weak-null sequence $(x_n) \subset E$ so that 
$\Vert x_n\Vert \geq c > 0$ for $n \in {\bf N}$.  Then
$\Vert Vx_n\Vert \to 0$ as $n \to \infty$ by the Dunford-Pettis property of
$E$ for any  $V \in W(E)$.  Hence
$$
\Vert x_n - Vx_n\Vert \geq c - \Vert Vx_n\Vert \geq {c\over 2}
$$
for all large enough $n \in {\bf N}$, so that $E$ does not have the W.A.P.
In particular, if $E$ contains an infinite-dimensional reflexive subspace,
then $E$ cannot have the Schur property.

\smallskip

\noindent (ii) Let $Q: E \to E/M$ stand for the weakly
compact quotient map. 
Suppose that there is a sequence
$(V_n) \subset W(E)$ so that
$\Vert Q - QV_n \Vert  \to 0$ as $n \to \infty$.
It follows that $QV_n$ is a compact operator $E \to E/M$ for 
$nÊ\in {\bf N}$, since $E$ has the DPP.
Hence the quotient map
$Q$ is a compact operator onto $E/M$, which is not  possible. 
$\square$

\medskip

Note that $c_0$ has the property that
every infinite-dimensional subspace $M \subset c_0$
fails to have the W.A.P. This follows from Proposition 6.10.(i) 
and the fact that $c_0$ is complementedly $c_0$-saturated, 
see [LT,2.a.2]. This fact is another point of difference between the W.A.P.
and classical approximation properties.

It is clear that $E$ has the W.A.P. 
if  $E$ has the Schur property and the B.A.P., since $W(E) = K(E)$
in this case.  Any space $E$ with the Schur property is $\ell^1$-saturated 
by Rosenthal's $\ell^1$-theorem (see [LT,2.e.5]).
This fact suggests the following question. 

\medskip

\noindent {\sl Question 6.11.\quad  Suppose that $E$ is an $\ell^1$-saturated Banach 
space that has the B.A.P.  Does $E$ have the W.A.P.?
}

\medskip

We answer Question 6.11 by showing that the Lorentz 
sequence spaces $d(w,1)$ fail to have the W.A.P.
Let $w = (w_j)$ be a positive non-increasing sequence satisfying
$$
w_1 =1, \quad \lim_{j\to\infty} w_j = 0 \quad {\rm and}\  \sum_{j=1}^{\infty} w_j = \infty . 
\eqno{(6.15)}
$$
Recall that $d(w,1)$  consists of the scalar sequences
$x = (x_j)$ for which
$$
\Vert x \Vert = \sum_{j=1}^{\infty} w_jx_j^* < \infty , 
$$
where $(x_j^*)$ is the non-increasing rearrangement of $(\vert x_j\vert )$.
The space $d(w,1)$ is $\ell^1$- saturated by [LT,4.e.3],  
but $d(w,1)$ does not have the DPP, since the coordinate basis 
$(e_n)$ and its biorthogonal sequence $(e_n^*)$ in $d(w,1)^*$ are weakly null.
In place of Proposition 6.10 we will use  
the (sub)symmetry of the Schauder basis $(e_n)$ for $d(w,1)$.

Let $E$ be a Banach space. Recall that a Schauder basis $(e_n)$ for $E$
is {\it symmetric} if $(e_{\pi (n)})$ and $(e_n)$ are equivalent for all 
permutations $\pi$ of ${\bf N}$. The basis $(e_n)$
is  {\it subsymmetric}, if $(e_n)$
is unconditional and $(e_{m_{n}})$ is equivalent to
$(e_n)$ for all subsequences $m_1 < m_2 < \ldots$. Every symmetric basis 
is also subsymmetric [LT,3.a.3].
Let $(x_j)$ and $(y_j)$ be basic sequences in  $E$.
Recall that  $(x_j)$ is said to {\it dominate} $(y_j)$ if 
$\sum_{j=1}^{\infty} c_jy_j$ converges in $E$ whenever
$\sum_{j=1}^{\infty} c_jx_j$ converges in $E$.

\medskip

\noindent {\it Example 6.12.}\quad  {\sl $d(w,1)$ does not have the W.A.P.
}

\medskip

\noindent {\it Proof.}  
The set $D = \{e_n: n \in {\bf N}\} \cup \{0\} \subset  d(w,1)$ 
is  weakly compact, since $(e_n)$ is a weak-null sequence. 
We will show that $D$ cannot be approximated in the sense of (1.1).
Suppose for this purpose that $V \in L(d(w,1))$ satisfies 
$\sup_{n \in {\bf N}} \Vert e_n - Ve_n\Vert
< {1\over {10}}$, and put $x_n = Ve_n$ for $n \in {\bf N}$.
Then the sequence $(x_n)$ is semi-normalized and weakly null.

\medskip

\noindent {\it Claim.}\quad $V \notin W(d(w,1))$.  

\medskip

First choose a basic subsequence 
$(x_{n_{j}})$ so that  $(x_{n_{j}})$ is equivalent to a block basic sequence $(y_j)$
of $(e_n)$, where $\Vert x_{n_{j}} - y_j\Vert \to 0$ as $j \to \infty$. 
Put $y_j = \sum_{k=p_{j}}^{q_j} a_ke_k$ for  $j \in {\bf N}$,
where $p_1 < q_1 < p_2 < q_2 < \ldots$ is a suitable sequence.
It is obvious that
$(e_{n_{j}})$ dominates $(x_{n_{j}}) = (Ve_{n_{j}})$. 
We next verify that 
$(x_{n_{j}})$ dominates $(e_{n_{j}})$, so that 
$(x_{n_{j}})$ and $(e_{n_{j}})$ will be equivalent
basic sequences in $d(w,1)$.  

We may assume by approximation that 
$n_j \in [p_j,q_j]$ and 
$a_{n_{j}} = e_{n_{j}}^*(y_j) > {9\over {10}}$ for $j \in {\bf N}$.
Let  $c_1,\ldots ,c_r$ be scalars and $r \in {\bf N}$.
It follows from the $1$-unconditionality of the basis $(e_j)$ that
$$
\Vert \sum_{j=1}^r c_jy_j\Vert = \Vert  \sum_{j=1}^r (\sum_{k=p_{j}}^{q_j} c_ja_ke_k)\Vert
\geq {9 \over {10}} \Vert \sum_{j=1}^r c_je_{n_{j}}\Vert ,
$$
and hence $(x_{n_{j}})$ dominates $(e_{n_{j}})$.

It follows that the restriction of $V$ determines a linear isomorphism
$[e_{n_{j}}: j \in {\bf N}] \to [Ve_{n_{j}}: j \in {\bf N}]$, since 
the sequences $(e_{n_{j}})$ and $(Ve_{n_{j}})$ are equivalent.
Here $[e_{n_{j}}: j \in {\bf N}] \approx d(w,1)$, because $(e_n)$
is a (sub)symmetric basis. This implies the Claim..
$\square$

\medskip

\noindent {\sl Remark 6.13.}\  The argument 
of Example 6.12 actually yields a more general observation,
which applies e.g. to certain Orlicz 
sequence spaces (see Chapter 4 of [LT]):

\smallskip 

{\sl Suppose that $E$ is a non-reflexive Banach space 
which has a weak-null, subsymmetric Schauder basis $(e_n)$.
Then $E$ does not have the W.A.P.}

\medskip
  
Azimi and Hagler  [AH]  introduced a class of spaces 
that provides a second solution to Question 6.11
(with some additional properties).
Let $w = (w_j)$ be a positive non-increasing sequence satisfying (6.15).
The Azimi-Hagler space $X(w)$ consists of the scalar sequences
$x = (x_j)$ for which
$$
\Vert x \Vert = 
\sup_{n; F_{1} < \ldots < F_{n}} \sum_{j=1}^n w_j \vert \sum_{k\in F_{j}} x_k\vert
 < \infty .  \eqno{(6.16)}
$$
The supremum is taken over all finite intervals $F_1 < \ldots < F_n$ of ${\bf N}$
and $n \in {\bf N}$.
The Banach space $X(w)$ is $\ell^1$- saturated, but it does not have 
the Schur property, see [AH,Thm. 1].
One point of interest in $X(w)$ comes from the facts that
the coordinate basis $(e_n)$ is not a subsymmetric basis for $X(w)$ 
(see the Remark on [AH,p. 295]), 
and $(e_n)$ does not even contain any weakly convergent subsequences 
(see [AH,Thm. 1.(3)]).
Hence the approach of Example 6.12 must be refined.

Let $P_{m,n}$ denote the natural projection of $X(w)$ onto $[e_s: m \leq s \leq n]$
for $m \leq n$. Thus $\Vert P_{m,n}\Vert \leq 2$.

\medskip

\noindent {\it Example 6.14.}\quad  {\sl $X(w)$ does not have the W.A.P.
}

\medskip

\noindent {\it Proof.} 
Put $z_n= e_{2n} - e_{2n-1}$ for $n \in {\bf N}$.
Then $(z_n)$ is a weak-null sequence in $X(w)$ (see [AH,Lemma 6]),
so that $\{z_n: n \in {\bf N}\} \cup \{0\}$ is  a weakly compact set.
Suppose that $V \in L(X(w))$ satisfies
$$
\sup_{n \in {\bf N}} \Vert z_n - Vz_n\Vert < {1\over {10}},
$$
and  set $x_n = Vz_n$ for $n \in {\bf N}$.
Thus $(x_n)$ is a semi-normalized weak-null sequence.

\medskip

\noindent {\it Claim.} {\sl $V \notin  W(X(w))$.}

\medskip

By the standard gliding hump argument 
we may first choose a subsequence $(x_{n_{j}})$
and natural numbers $p_1 < q_1 < p_2 < q_2 < p_3 < \ldots$, so that

\smallskip

\item{(i)} $(x_{n_{j}})$ and $(y_j)$ are equivalent basic sequences,
where $y_j = P_{p_{j},q_{j}}(x_{n_{j}})$ for $j \in {\bf N}$,

\smallskip

\item{(ii)}  $y_j = u_j + a_je_{2n_{j}-1} + b_je_{2n_{j}} + v_j$, where
$p_j < 2n_{j}-1 < 2n_j < q_j$, and the supports satisfy
$supp(u_j) \subset [p_j,2n_{j}-1)$ and $supp(v_j) \subset (2n_j,q_j]$
for $j \in {\bf N}$.

\smallskip

\item{(iii)}  $\vert a_j +1\vert < {2 \over {10}}$ and  $\vert b_j -1\vert < {2 \over {10}}$
for $j \in {\bf N}$.

\smallskip

\noindent Property (iii) follows from the fact that $\Vert y_j - z_{n_{j}}\Vert  = 
\Vert P_{p_{j},q_{j}}(x_{n_{j}} - z_{n_{j}})\Vert < {2 \over {10}}$.

Clearly $(z_{n_{j}})$ dominates $(x_{n_{j}}) = (Vz_{n_{j}})$.
By property (i) it suffices to verify that there is $c > 0$ so that
$$
\Vert \sum_j c_jy_j\Vert \geq c \cdot \Vert \sum_{j} c_jz_{n_{j}}\Vert \eqno{(6.17)}
$$
for all scalars $c_1, c_2, \ldots$. Indeed, in that event  
the basic sequences $(x_{n_{j}})$ 
and $(z_{n_{j}})$ are equivalent, and 
the fact that $X(w)$ is $\ell^1$-saturated [AH, Thm. 1] 
will imply that $V$ fixes some $\ell^1$-copy 
contained in $[z_{n_{j}}: j \in {\bf N}]$.

It is enough to verify (6.17)
for all finite sums $z = \sum_{j=1}^r c_jz_{n_{j}}$ and $y =  \sum_{j=1}^r c_jy_j$.
Put $F^*(x) = \sum_{s\in F} x(s)$ for $x = \sum_{s=1}^{\infty} x(s)e_s \in X(w)$,
whenever $F \subset {\bf N}$ is a finite interval.
Suppose that $F_1 < F_2 < \ldots < F_m$ are finite intervals
for which $\sum_{i=1}^m w_i\vert F_i^*(z)\vert = \Vert z \Vert$.
The non-zero terms $\vert F_i^*(z)\vert$ have the form
$\vert c_j\vert$ or $\vert c_j - c_k\vert$ for suitable $j < k$. Indeed, 
there is no contribution to $F_i^*(z)$ from the terms $c_l(e_{2n_{l}-1} - e_{2n_{l}})$, 
where both  $2n_{l}-1, 2n_l \in F_i$. 

If $i$ is such that $\vert F_i^*(z)\vert = \vert c_j\vert$,
then we may replace $F_i$ by $G_i = \{2n_{j} - 1\}$ or $G_i = \{2n_{j}\}$ 
without affecting $\vert F_i^*(z)\vert = \vert c_j\vert$.
The choice of $2n_{j} - 1$ or $2n_{j}$ is according to which of these indices 
contributes the term $\vert c_j\vert$.
Thus $G_i^*(y) = b_jc_j$ or $G_i^*(y) = a_jc_j$, so that (iii) yields
$$
\vert G_i^*(y)\vert \geq \min\{\vert a_j\vert , \vert b_j\vert \} \cdot \vert c_j\vert 
\geq {8\over {10}} \vert c_j\vert = {8\over {10}} \vert F_i^*(z)\vert . \eqno{(6.18)}
$$
If $\vert F_i^*(z)\vert = \vert c_j - c_k\vert$ for some $j < k$, then we
consider  two singletons $G_{i,1} < G_{i,2}$ of the preceding type instead of
$F_i$. In this case we get as above that
$$
\vert F_i^*(z)\vert = \vert c_j - c_k\vert \leq 
\vert c_j\vert + \vert c_k\vert \leq {{10}\over 8} (\vert G_{i,1}^*(y)\vert 
+ \vert G_{i,2}^*(y)\vert ). \eqno{(6.19)}
$$
Put $A = \{ i: \vert F_i^*(z)\vert = \vert c_j\vert$ for some $j\}$ and
$B = \{ i: \vert F_i^*(z)\vert = \vert c_j - c_k\vert$ for some $j < k\}$.
We get from (6.18) and (6.19) that
$$
\eqalign{\Vert z \Vert = \ & 
\sum_{i\in A} w_i\vert F_i^*(z)\vert + \sum_{i\in B} w_i\vert F_i^*(z)\vert \cr
\leq \ & {{10}\over 8} (\sum_{i\in A} w_i\vert G_i^*(y)\vert + 
\sum_{i\in B} w_i\vert G_{i,1}^*(y)\vert )
+  {{10}\over 8} \sum_{i\in B} w_i\vert G_{i,2}^*(y)\vert 
\leq  {{20}\over 8}
\Vert y\Vert . \cr
}
$$
Here $\sum_{i\in B} w_i\vert G_{i,2}^*(y)\vert \leq \Vert y\Vert$,
since the weight sequence $(w_i)$ is non-increasing.
$\square$

\medskip

\noindent {\sl Remark 6.15.}\  The arguments for Examples 6.12 and 6.14 
demonstrate that
$d(w,1)$ and $X(w)$ even fail to have the unbounded W.A.P.

\medskip

The inner W.A.P.  (see section 5) is more difficult to study.
Our final example shows that the dual $JT^*$ of the James tree space
does not have the inner W.A.P.
We will require the following facts from [LS,Thm. 1]:
$JT$ has a predual $B$ and $B^{**}/B = \ell^2(\Gamma )$
isometrically, where 
$\Gamma$ is the uncountable collection of all branches of
${\cal T}$. In particular,
$JT^{***}/JT^* = \ell^2(\Gamma )$. 

\medskip

\noindent {\sl Example 6.16. \quad
$JT^*$ does not have the inner WAP.
}

\medskip

\noindent {\it Proof.} The argument is a modification of that of 
[T2,1.4] for the Johnson-Lindenstrauss space. 
Let  $q: JT^* = B^{**} \to \ell^2(\Gamma )$ be the weakly compact quotient map.  
Suppose to the contrary that there is a 
sequence $(V_n) \subset
W(JT^*)$ such that 
$$
\lim_{n\to \infty} \Vert q - qV_n\Vert  = 0. \eqno{(6.20)}
$$
Recall next that any weakly compact set $D \subset JT^*$
is norm separable, since $(D,w)$ is metrizable in this case.
This is based on the fact that
$JT$ is a separable space not containing any
copies of $\ell^1$ (see [LS,Cor. 1] or  [FG,3.a.8]), so that 
$JT$ is $w^*$-sequentially dense in $JT^{**}$ by the Main Theorem of [OR].

Deduce that the closure 
$\overline{q(V_nB_{JT^*})}$ is a norm separable set in $\ell^2(\Gamma )$
for $n \in {\bf N}$. Thus (6.20) implies that 
$B_{\ell^2(\Gamma )} = \overline{q(B_{JT^*})}$ is also norm separable
by approximation.
This contradicts the non-separability of 
$\ell^2(\Gamma )$.   $\square$

\medskip
 
If $E$ has the inner W.A.P., 
then $W(E)$ has a B.R.A.I. by [BD,Prop. 11.2]  (cf. also the proof of 
Proposition 2.5.(i)).  
Thus Example 6.16 and [LW,Cor. 2.4]  suggest

\medskip

\noindent {\sl Problems 6.17. (i) Does $J$ have the {\it  inner} W.A.P.?
(ii) Does $JT$ have the inner W.A.P.?
(iii) Does $JT^*$ have the W.A.P.?
}

\medskip

\noindent {\sl Acknowledgements.}  The second-named author is greatly 
indebted to E. Odell and H.P. Rosenthal for their hospitality
during his visit to the Department of Mathematics,
the University of Texas at Austin, in late 2002.

\bigskip

\noindent {\bf References.}

\medskip

\item{[A]}  A. Andrew: James' quasi-reflexive space is not
isomorphic to any subspace of its dual, {\it Israel J. Math.} 38
(1981), 276-282.

\item{[ArT]}  S.A. Argyros and A. Tolias:  Methods in the 
theory of hereditarily indecomposable Banach spaces,
{\it Mem.  Amer. Math. Soc.} (to appear).

\item{[AT]} K. Astala and H.-O. Tylli: Seminorms related 
to weak compactness and to Tauberian operators,
{\it Math. Proc. Cambridge Phil. Soc.} 107 (1990), 367-375.

\item{[AH]} P. Azimi and J.N. Hagler: Examples of hereditarily
$\ell^1$ Banach spaces failing the Schur property, {\it Pacific J. Math.}
122 (1986), 287-297.

\item{[BHO]} S.F. Bellenot, R. Haydon and E. Odell: Quasi-reflexive and 
tree spaces constructed in the spirit of R.C. James, {\it Contemp. Math.}
85 (1989), 19-43.

\item{[BD]} F.F. Bonsall and J. Duncan: {\it Complete normed algebras.}
Ergebnisse der Mathematik vol. 80 (Springer, 1973).

\item{[B]} J. Bourgain: {\it New classes of ${\cal L}^p$-spaces.}
Lecture Notes in Mathematics vol. 889 (Springer-Verlag, 1981).

\item{[CLL]} P.G. Casazza, B.L. Lin and R.H. Lohman: On James'
quasi-reflexive Banach space, {\it Proc. Amer. Math. Soc. } 67 (1977),
265-271.

\item{[D]} H.G. Dales: {\it Banach algebras and automatic continuity.}
London Mathematical Society Monographs vol. 24 (Oxford University Press, 2000).

\item{[Di]} J. Diestel:  A survey of results related to the Dunford-Pettis property,
{\it Contemp. Math.} 2 (1980), 15-60.

\item{[FG]} H. Fetter and B. Gamboa de Buen: {\it The James
Forest}, London Math. Soc. Lecture Notes 236  (Cambridge University Press, 1997).

\item{[GST]} M. Gonzalez, E. Saksman and H.-O. Tylli:
Representing non-weakly compact operators,
{\it Studia Math.} 113 (1995), 265-282.

\item{[GW]} N. Gr{\o}nb{\ae}k and G.A. Willis: Approximate identities in 
Banach algebras of compact operators. {\it Canad. Math. Bull.} 36 (1993), 45-53.

\item{[J1]}  R.C. James: {\it A non-reflexive Banach space isometric with
its second conjugate space.} Proc. Nat. Acad. Sci. U.S.A. 37 (1951),
174-177.

\item{[J2]} R.C. James: A separable somewhat reflexive
Banach space with nonseparable dual, {\it Bull. Amer. Math. Soc.} 
80 (1974), 738-743.

\item{[L]} I.E. Leonard: Banach sequence spaces, {\it J. Math. Anal. Appl.}
54 (1976), 245-265.

\item{[LNO]} {\AA}. Lima, O. Nygaard and E. Oja: Isometric factorization
of weakly compact operators and the approximation property,
{\it Israel J. Math.} 119 (2002), 325-348.

\item{[LS]} J. Lindenstrauss and C. Stegall: Examples of separable spaces
which do not contain $\ell^1$ and whose duals are not separable, 
{\it Studia Math.} 54 (1975), 81-105.

\item{[LT]} J. Lindenstrauss and L. Tzafriri: {\it Classical Banach spaces I.
Sequence spaces.} Ergebnisse der Mathematik vol. 92 (Springer, 1977).

\item{[LW]} R.J. Loy and G.A. Willis: Continuity of derivations
on $B(E)$ for certain Banach spaces $E$, {\it J. London Math. Soc.} 40 (1989),
327-346.

\item{[OR]}  E. Odell and H.P. Rosenthal:  A double dual characterization of
separable Banach spaces containing $\ell^1$, {\it Israel J. Math.}
20 (1975), 375-384.

\item{[Ph]}  R.R. Phelps: {\it Lectures on Choquet's Theorem} (2. edition).
Lecture Notes in Mathematics vol. 1757 (Springer, 2001).

\item{[P]} G. Pisier: The dual $J^*$ of the James space has cotype $2$ 
and the Gordon Lewis property, {\it Math. Proc. Cambridge Phil.  Soc.} 103 (1988),
323-331.

\item{[PQ]} G. Pisier and Q. Xu:  Random series in the interpolation
spaces between the spaces $v_p$, in {\it Geometric Aspects of Functional 
Analysis 1985-86} (J. Lindenstrauss and V. Milman, eds.) 
Lecture Notes in Mathematics vol. 1267 (Springer, 1987),
pp. 185-209.

\item{[Pt]} V.  Ptak: A combinatorial theorem on systems of inequalities
and its applications to analysis, {\it Czech. Math. J.} 84 (1959), 629-630.

\item{[R]} O. Reinov: How bad can a Banach space with the approximation
property be ? {\it Math. Notes} 33 (1983),  427-434.

\item{[T1]} H.-O. Tylli: The essential norm of an operator
is not self-dual, {\it Israel J. Math.} 91 (1995), 93-110.

\item{[T2]} H.-O. Tylli: Duality of the weak essential norm,
{\it Proc. Amer. Math. Soc.} 129 (2001), 1437-1443.

\item{[W]} M. Wojtowicz: On the James space $J(X)$ for a Banach space $X$,
{\it Comment. Math. Prace Mat.} 23 (1983), 183-188.

\bigskip

\noindent Addresses:

\noindent (Odell)

\noindent Department of Mathematics

\noindent The University of Texas at Austin

\noindent Austin, TX 78712

\noindent USA

\smallskip

\noindent e-mail: {\it odell@math.utexas.edu}

\medskip

\noindent (Tylli)

\noindent Department of Mathematics

\noindent P.B. 4 (Yliopistonkatu 5)

\noindent FIN-00014 University of Helsinki

\noindent Finland

\smallskip

\noindent e-mail: {\it hojtylli@cc.helsinki.fi}

\end